\documentclass[11pt]{article}
\usepackage{amsfonts}
\usepackage{color}
\newtheorem{theo}{Theorem}[section]
\newtheorem{lem}[theo]{Lemma}
\newtheorem{cor}[theo]{Corollary}

\newtheorem{prop}[theo]{Proposition}

\newcommand{\mysection}[1]{\section{#1} \setcounter{equation}{0}}
\newcommand{\proof}{{\sc Proof.} \quad}
\newcommand{\proofc}{{\sc Proof} \ }
\newcommand{\be}{\begin{equation} \label}
\newcommand{\ee}{\end{equation}}
\newcommand{\bea}{\begin{eqnarray}\label}
\newcommand{\eea}{\end{eqnarray}}
\newcommand{\bas}{\begin{eqnarray*}}
\newcommand{\eas}{\end{eqnarray*}}
\newcommand{\bit}{\begin{itemize}}
\newcommand{\eit}{\end{itemize}}
\newcommand{\qed}{\hfill$\Box$ \vskip.2cm}
\newcommand{\nn}{\nonumber}
\newcommand{\R}{\mathbb{R}}
\newcommand{\N}{\mathbb{N}}

\newcommand{\eps}{\varepsilon}

\newcommand{\supp}{{\rm supp} \, }

\newcommand{\abs}{\\[5pt]}
\newcommand{\Abs}{\\[5mm]}

\newcommand{\mint}{-\hspace*{-4mm} \int}
\newcommand{\tm}{T_{max}}
\newcommand{\io}{\int_\Omega}
\newcommand{\tu}{\widetilde{u}}
\newcommand{\tU}{\widetilde{U}}
\newcommand{\UU}{\underline{U}}
\newcommand{\OU}{\overline{U}}
\newcommand{\kz}{K_0}
\newcommand{\parab}{{\cal P}}
\newcommand{\qarab}{{\cal Q}}
\begin{document}
\title{Critical mass for infinite-time aggregation in a chemotaxis model with indirect signal production}
\author{
Youshan Tao\footnote{taoys@dhu.edu.cn}\\
{\small Department of Applied Mathematics, Dong Hua University,}\\
{\small Shanghai 200051, P.R.~China}
\and
Michael Winkler\footnote{michael.winkler@math.uni-paderborn.de \ (corresponding author)}\\
{\small Institut f\"ur Mathematik, Universit\"at Paderborn,}\\
{\small 33098 Paderborn, Germany} }
\date{}
\maketitle
\begin{abstract}
\noindent
  We study the Neumann initial-boundary problem for the chemotaxis system
  \bas
    \left\{ \begin{array}{ll}
    u_t= \Delta u - \nabla \cdot (u\nabla v), & x\in \Omega, \, t>0, \\[1mm]
    0=\Delta v - \mu(t)+w,
    & x\in \Omega, \, t>0, \\[1mm]
    \tau w_t + \delta w = u,
    & x\in \Omega, \, t>0,
    \end{array} \right.
    \qquad \qquad (\star)
  \eas
  in the unit disk $\Omega:=B_1(0)\subset \R^2$, where $\delta\ge 0$ and $\tau>0$ are given parameters and
  $\mu(t):=\mint_\Omega w(x,t)dx$, $t>0$.\abs
  It is shown that this problem exhibits a novel type of critical mass phenomenon with regard to the formation
  of singularities, which drastically differs from the well-known threshold property of the classical Keller-Segel system,
  as obtained upon formally taking $\tau\to 0$,
  in that it refers to blow-up in infinite time rather than in finite time:\\
  Specifically, it is first proved that for any sufficiently regular nonnegative initial data $u_0$ and $w_0$, ($\star$)
  possesses a unique global classical solution. In particular, this shows that in sharp contrast
  to classical Keller-Segel-type systems reflecting immediate signal secretion by the cells themselves,
  the indirect mechanism of signal production in ($\star$) entirely rules out any occurrence of blow-up in finite time.\\
  However, within the framework of radially symmetric solutions it is next proved that
  \bit
  \item
    whenever $\delta>0$ and $\io u_0<8\pi\delta$, the solution remains uniformly bounded, whereas
  \item
    for any choice of $\delta\ge 0$ and $m>8\pi\delta$, one can find initial data such that $\io u_0=m$, and such that
    for the corresponding solution we have
    \bas
    \|u(\cdot,t)\|_{L^\infty(\Omega)} \to \infty
    \qquad \mbox{as } t\to\infty.
    \eas
  \eit
 {\bf Key words:} chemotaxis, infinite-time blow-up, critical mass\\
 {\bf AMS Classification:} 35B44, 35B51 (primary); 35A01, 35K55, 35Q92, 92C17 (secondary)
\end{abstract}
\newpage
\mysection{Introduction}
{\bf A chemotaxis model with indirect signal production.} \quad
Chemotaxis, the biased movement of cells along
concentration gradients of a chemical signal, is known to play a significant role
in numerous biological circumstances such as
bacterial aggregation, spatial pattern formation, embryonic
morphogenesis, cell sorting, immune response, wounding healing,
tumor-induced angiogenesis, and also tumor invasion
(see \cite{tindall2008}, \cite{maini_1991}, \cite{painter_bmb2009},
\cite{friedman_hu_xue}, \cite{anderson_chaplain},
\cite{fontelos_friedman_hu}, \cite{chaplain_lolas_2005} and
\cite{chaplain_lolas_2006}, for instance). The renowned Keller-Segel
model (cf.~(\ref{2}) below), describing the collective behavior of cells in response to a
signal produced by the cells themselves, has been well-studied
with regard to biological implications, but beyond this, during the last decades
quite a thorough comprehension of its mathematical features has grown
in various directions (\cite{tindall2008}, \cite{hillen_painter2009}, \cite{BBTW_survey}).\abs
In contrast to this well-understood paradigmatic case, the theoretical understanding
is much less developed in situations when a chemotactic cue is not released by the cells themselves.
Typical examples for such mechanisms include cases when the signal is not produced at all, such as in
oxygenotaxis processes of swimming aerobic bacteria
which preferably move toward higher concentrations of externally provided oxygen as their nutrient
(\cite{goldstein}),
and also cases in which signal production occurs within more complex processes, possibly involving chemical
reactions or even cascades thereof, such as e.g.~in the glycolysis reaction
(\cite{dillon_maini_othmer}, \cite{painter_maini_othmer};
cf.~also \cite{murray} and \cite{calvez_perthame} for further extensions of chemotaxis models involving
additional couplings.\abs
It is the purpose of the present work to
achieve some insight into possible features of chemotaxis models
accounting for the latter type of more complex signal production mechanisms.
Specifically, we shall be concerned with the prototypical parabolic-elliptic-ODE system
\be{0}
    \left\{ \begin{array}{ll}
    u_t= \Delta u - \nabla \cdot (u\nabla v), & x\in \Omega, \, t>0, \\[1mm]
    0=\Delta v - \mu(t)+w,
    & x\in \Omega, \, t>0, \\[1mm]
    \tau w_t + \delta w = u,
    & x\in \Omega, \, t>0, \\[1mm]
    \frac{\partial u}{\partial\nu}=\frac{\partial v}{\partial\nu}=0,
    & x\in \partial\Omega, \, t>0, \\[1mm]
    u(x,0)=u_0(x), \quad w(x,0)=w_0(x),
    & x\in \Omega,
    \end{array} \right.
\ee
in the unit disk $\Omega:=B_1:=B_1(0)\subset \R^2$,
where $\delta\ge 0$ and $\tau>0$ are fixed parameters and\be{mu_def}
    \mu(t):=\mint_\Omega w(x,t)dx, \qquad t>0.
\ee
In a concrete biological framework,
this model arises as a simplification of the chemotaxis model recently proposed
by Strohm, Tyson and Powell in \cite{powell} to describe the spread
and aggregative behavior of the Mountain Pine Beetle (MPB) in a forest
habitat
considered negligibly thin in its vertical dimension.
Their model involves three variables: the density of flying
MPB, denoted by $u$, the density of nesting MPB, represented by $w$,
and the concentration $v$ of beetle pheromone, the latter being secreted only
by those MBP which are nested in trees.
Besides random diffusive motion, the flying MPB can partially orient their
movement according to concentration gradients of MPB pheromone. Once MPB
nest they do not move any longer, thus meaning that apart from the
increase of $w$ through transition from the flying to the nested state,
the only further quantity relevant to their evolution remains their death rate $\delta$.
For more details on the physical background, we refer the reader to \cite[Section 2]{powell}.\abs
From a mathematical point of view, (\ref{0}) can be viewed as a variant of the Keller-Segel model
associated with the system
 \be{1}
    \left\{ \begin{array}{ll}
    u_t= \Delta u - \nabla \cdot (u\nabla v), & x\in \Omega, \, t>0, \\[1mm]
    0=\Delta v - \widetilde{\mu} +u,
    & x\in \Omega, \, t>0,
    \end{array} \right.
\ee where $\widetilde{\mu} :=\mint_\Omega u \equiv \mint_\Omega
u_0$. which can formally be obtained from (\ref{0}) upon taking
$\tau\searrow 0$. In the case when $\Omega$ coincides with the
entire space $\R^2$ or $\R^3$, the corresponding limit system of the
latter arises in the modeling of self-gravitating particles
(\cite{biler_hilhorst_nadzieja}), and furthermore it was introduced
in \cite{JL} as a simplification of the well-known classical
Keller-Segel model (\cite{KS}) of chemotaxis, the original version
of which being \be{2}
    \left\{ \begin{array}{ll}
    u_t= \Delta u - \nabla \cdot (u\nabla v), & x\in \Omega, \, t>0, \\[1mm]
    v_t=\Delta v - v+u,
    & x\in \Omega, \, t>0.
    \end{array} \right.
\ee
Here the hypotheses justifying the reduction of (\ref{2}) to (\ref{1}), namely the
physically meaningful assumptions that chemicals diffuse much faster than cells, and that the
particular signal substance in question degrades sufficiently slowly, have been used in various related
contexts and are also part of the simplification of the original model in \cite{powell} to (\ref{0})
(cf.~also the review paper \cite{horstmann_dmv}).\abs
Let us emphasize here the evident difference between (\ref{0}) for $\tau>0$ on the one hand and the two-component
Keller-Segel systems (\ref{1}) and (\ref{2}) on the other:
In both of the latter, the quantity $u$ {\em directly} produces the quantity $v$ governing its cross-diffusion, 
whereas the corresponding signal production in (\ref{0})
occurs in an {\em indirect} process, with first $u$ producing the
third quantity $w$, and with the latter being exclusively
responsible for the release of $v$.\Abs
{\bf Blow-up and critical mass phenomena.} \quad
It is known that chemotactic cross-diffusion terms, constituting
the apparently most characteristic model ingredient in all systems (\ref{0}), (\ref{1}) and (\ref{2}),
may have a strong destabilizing potential and even enforce the formation of singularities.
Correspon\-dingly, a striking feature of both Keller-Segel systems (\ref{1}) and (\ref{2})
appears to be the occurrence of some solutions blowing up in finite time, which is commonly
viewed as mathematically expressing numerous processes of spontaneous cell aggregation which can be observed in
experiments (see \cite{hillen_painter2009} and also \cite{BBTW_survey} for a survey).
Indeed, in the spatially two-dimensional framework considered here, the appearance of
such explosion phenomena is closely related to the initially present total mass $\io u_0$ of cells.
For instance, it was shown in \cite{JL} and \cite{biler1998} that in
the spatially radial setting, the system (\ref{1}) possesses some
solutions which blow up in finite time provided that this mass $\io
u_0$ is large enough, whereas solutions remain bounded whenever $\io
u_0$ is small; as a precise value distinguishing the respective mass
regimes either allowing for or suppressing explosions the {\em
critical mass} $m_c=8\pi$ could be identified (cf.~\cite{biler1998},
\cite{nagai2001} and \cite{perthame} for (\ref{1}) and closely
related variants thereof).\abs
As for the fully parabolic chemotaxis system (\ref{2}), an analogous critical mass phenomenon is known to occur,
the respective threshold value again being $m_c=8\pi$ in the radially symmetric situation.
For corresponding results on boundedness in the subcritical regime we refer to \cite{nagai_senba_yoshida};
some quite particular blow-up solutions with $\io u_0>8\pi$ have been detected in \cite{herrero_velazquez},
whereas recently in \cite{mizoguchi_winkler} it was shown that such a singularity formation indeed occurs
within a considerably large set of supercritical-mass initial data, which can even be viewed generic in an
appropriate sense.\abs
In the nonradial setting, corresponding critical mass phenomena seem to be present, with
a reduced value of $m_c=4\pi$. For parabolic-elliptic Keller-Segel systems, rigorous proofs for this can be found in,
or easily adapted from \cite{nagai2001} and \cite{nagai_senba_yoshida}; in the parabolic-parabolic case,
only a respective boundedness result is available in the case $\io u_0<4\pi$ (\cite{nagai_senba_yoshida}),
whereas the occurrence of any nonradial finite-time blow-up solution to (\ref{2}) appears to be a challenging
open problem (cf.~\cite{horstmann_wang} for a partial result on unboundedness).\abs
Let us mention that in the spatially one-dimensional versions of
both (\ref{1}) and (\ref{2}), all solutions emanating from
conveniently smooth initial data are global in time and remain
uniformly bounded (\cite{osaki_yagi}), while in three- or
higher-dimensional balls, for arbitrarily small values of $m>0$ one
can find smooth initial data fulfilling $\io u_0=m$, for which the
corresponding solution will blow up in finite time (see
\cite{nagai_senba1996} for a parabolic-elliptic and \cite{win_jmpa}
for the fully parabolic case). A critical mass phenomenon thus
occurs only in the two-dimensional situation.\Abs
{\bf Main results. A novel type of critical mass phenomenon.} \quad
It is the purpose of the present paper to rigorously
investigate the qualitative features of the system (\ref{0}) with regard to its original intention to model
processes of aggregation.
Here our focus will be on the question in how far the indirect signal production mechanism in (\ref{0})
can enforce singularity formation in the first solution component $u$.
Our main results in this direction
show that actually also (\ref{0}) exhibits a type of critical mass phenomenon, but that the latter appears to be
novel in the context of chemotaxis problems: Surprisingly, namely,
unlike that for (\ref{1}) and (\ref{2}), the mass threshold property we shall identify here will
refer to blow-up {\em in infinite time} rather than in finite time.\abs
Indeed, by deriving energy-type estimates through rather straightforward testing procedures
we can first show that for all reasonably regular initial data with arbitrary mass $\io u_0$,
(\ref{0}) is globally classically solvable:
\begin{prop}\label{prop100}
  Let $\delta\ge 0$ and $\tau>0$, and
  suppose that $u_0\in C^0(\bar{\Omega})$ and $w_0\in C^1(\overline{\Omega})$   
  are nonnegative.
  Then there exists a unique triple $(u, v, w)$ of nonnegative functions
  \bas
    & & u\in C^0(\overline{\Omega}\times [0,\infty)) \cap C^{2,1}(\overline{\Omega}\times (0,\infty)), \\
    & & v\in C^{2,0}(\overline{\Omega}\times [0,\infty)), \\
    & & w \in C^{0,1}(\overline{\Omega}\times [0,\infty)),
  \eas
  which solves (\ref{0}) in the classical sense.
\end{prop}
We shall next detect the number
\bas
    m_c:=8\pi\delta
\eas
to be critical with regard to boundedness of radial solutions.
The first part of this characterization is contained in the following.
\begin{theo}\label{theo200}
  Let $\delta>0$ and $\tau>0$, and suppose that $u_0\in C^0(\overline{\Omega})$ and $w_0\in C^1(\overline{\Omega})$
  are radially symmetric and nonnegative, and that $m:=\io u_0$ satisfies
  \bas
    m<8\pi\delta.
  \eas
  Then the solution of (\ref{0}) is bounded in $\Omega\times (0,\infty)$; that is, there exists a constant $C>0$ such that
  \bas
    u(x,t) \le C, \quad v(x,t) \le C \quad \mbox{and} \quad w(x,t) \le C
    \qquad \mbox{for all $x\in\Omega$ and } t>0.
  \eas
\end{theo}
Secondly, the above picture is completed by our final statement: In fact, for any $m>8\pi\delta$
we shall derive an essentially explicit condition
on the radial initial data $u_0$ and $w_0$ which under the assumption $\io u_0=m$ ensures that in the large time limit,
the solution diverges exponentially in its first component when measured in $L^\infty(\Omega)$:
\begin{theo}\label{theo300}
  Let $\delta\ge 0$ and $\tau>0$. Then for any
  \bas
    m>8\pi\delta
  \eas
  there exist $R\in (0,1)$ and $\alpha>0$ such that for each
  $\eta>0$ one can find positive constants $\Gamma_u(m,\eta), \gamma(m,\eta)$ and $\Gamma_w(m,\eta)$
  with the property
  that for all radially symmetric nonnegative functions $u_0\in C^0(\overline{\Omega})$ and
  $w_0\in C^1(\overline{\Omega})$
  satisfying
  \be{300.00}
    \io u_0=m>8\pi\delta
  \ee
  and
  \be{300.01}
    \mint_{B_r} u_0 \ge \Gamma_u(m,\eta)
    \qquad \mbox{for all } r\in (0,R)
  \ee
  and
  \be{300.02}
    \mint_{B_1 \setminus B_r} u_0 \le \gamma
    \qquad \mbox{for all } r\in (R,1)
  \ee
  as well as
  \be{300.03}
    \mint_{B_r} w_0 \ge \mint_{B_1} w_0 + \Gamma_w(m,\eta)
    \qquad \mbox{for all } r\in (0,R)
  \ee
  and
  \be{300.04}
    \mint_{B_1\setminus B_r} w_0 \le \mint_{B_1} w_0 -\eta
    \qquad \mbox{for all } r\in (R,1),
  \ee
  the corresponding solution $(u,v,w)$ of (\ref{0}) is unbounded in the sense that
  \bas
    \|u(\cdot,t)\|_{L^\infty(\Omega)} \ge c e^{\alpha t}
    \qquad \mbox{for all } t>0
  \eas
  with some $c=c(m, \eta, \delta, \tau, \|w_0\|_{L^1(\Omega)})>0$.
\end{theo}
As a particular consequence, this provides some quantitative information on the damping role of the death rate
$\delta$ in (\ref{0}).
For instance, it follows from Theorem \ref{theo200} that for any given initial data $(u_0,w_0)$ the associated
solution will remain bounded whenever $\delta>0$ is suitably large.
On the other hand, if $\delta$ vanishes then unbounded solutions can be found for arbitrarily small values of the
initial mass $\io u_0$.\abs
Moreover, the criticality of $m_c=8\pi\delta$, as thus detected to predict the possibility or impossibility
of aggregation phenomena in (\ref{0}) for positive values of $\tau$ and $\delta$, appears to be consistent with the above
mass threshold properties of (\ref{1}): Indeed, in the limit case $\tau=0$, in which in (\ref{0}) clearly any
initial condition on $w$ becomes obsolete, we will have $w\equiv \frac{1}{\delta} u$.
Hence, upon substituting $\tilde u:=\frac{1}{\delta} u$ we see that we may assume that $\delta=1$, and that
(\ref{0}) reduces to (\ref{1}), having critical mass $m_c=8\pi=8\pi\delta$;
the fact that $m_c$ is then related to finite-time blow-up, rather than to inifinite-time aggregation,
may be viewed as a consequence of the lacking relaxation mechanism reflected in the ODE for $w$ in (\ref{0})
when $\tau>0$.
In summary, varying $\tau$ over the interval $[0,\infty)$ does not change the {\em value} of the critical mass,
but it significantly affects its precise {\em role} when passing from positive $\tau$ to the case $\tau=0$.\Abs
{\bf Main ideas underlying our approach.} \quad
Let us briefly outline the methods we pursue in the derivation of Theorem \ref{theo200} and Theorem \ref{theo300}.
Our approach to both of these will be based on a transformation reducing (\ref{0}) to an initial-boundary value problem
for a scalar degenerate parabolic equation.
Though well-established in related contexts, this transformation results in an equation which,
unlike the corresponding situation in the standard
Keller-Segel system (\ref{1}) (\cite{JL}), now contains a nonlinear production term that is nonlocal in time.
More precisely, we shall see that the mass distribution function $U$ associated with a given radial solution
$u=u(r,t)$ of (\ref{0}), that is, the function defined by
\bas
    U(\xi,t):=\int_0^{\sqrt{\xi}} ru(r,t) dr, \qquad \xi\in [0,1], \ t\ge 0,
\eas
satisfies the single equation
\be{foo}
    U_t= 4\xi U_{\xi\xi}
    + \frac{2}{\tau} \Bigg\{ \int_0^t e^{-\frac{\delta}{\tau}(t-s)} \Big( U(\xi,s)-\frac{m}{2\pi} \xi\Big) \, ds \Bigg\}
    \cdot U_\xi
    + 2\Big(W_0(\xi)-\kz \xi\Big) \cdot e^{-\frac{\delta}{\tau} t} U_\xi,
\ee
for $\xi\in (0,1)$ and $t>0$, where $W_0(\xi):=\int_0^{\sqrt{\xi}} rw_0(r)dr$, $\xi\in [0,1]$, and $\kz:=W_0(1)$
(cf. Lemma \ref{lem31}). Clearly, $u$ is bounded if and only if the spatial gradient $U_\xi$ is bounded.
Fortunately, the corresponding parabolic operator allows for a comparison principle
(Lemma \ref{lem12}), and thus enables us to focus our subsequent analysis on the construction
of appropriate super- and subsolutions.\abs
Based on such a comparison argument, under the subcriticality assumption $m<8\pi\delta$ from Theorem \ref{theo200}
we shall first obtain an estimate of the form $U(\xi,t) \le C\xi$ for all $(\xi,t)\in (0,1)\times (0,\infty)$ and some $C>0$
(Lemma \ref{lem10}).
This means that given $\eps>0$, adjusting $r_0\in (0,1)$ suitably we can achieve that
the mass which the original solution accumulates in the ball $B_{r_0}(0)$ satisfies $\int_{B_{r_0}(0)} u(x,t)dx<\eps$
for all $t>0$. In conjunction with a corresponding $\eps$-regularity result (Section \ref{sect5.3}) this will yield the
desired boundedness property of such solutions.\abs
In the case $m>8\pi\delta$ addressed in Theorem \ref{theo300}, we will construct subsolutions exhibiting
gradient grow-up at the origin; that is, we shall find a family of adequate
subsolutions $\UU$ to (\ref{foo}) with the properties $\UU(0,t)=0$ for all $t>0$ and
$\UU_\xi(0,t) \to + \infty$ as $t\to\infty$.
Proving Theorem \ref{theo300} then amounts to finding sufficient conditions for
$u_0$ and $w_0$ ensuring that $U(\xi,0)\ge \UU(\xi,0)$ for all $\xi\in (0,1)$.\\
We find it worthwhile to underline here that the structure near the origin of the latter
comparison functions, to be explicitly constructed and analyzed in detail in Section \ref{sect6}, will be given by
\be{foo1}
    \UU(\xi,t):=\frac{a(t)\xi}{b(t)+\xi}, \qquad \xi\in [0,\xi_0), \ t\ge 0,
\ee
with $b(t)=b_0 e^{-\alpha t}$, $t\ge 0$, and appropriately chosen $a\in C^1([0,\infty))$, $\xi_0\in (0,1)$, $b_0>0$
and $\alpha>0$.
The idea for this construction originates from standard knowledge on equilibria for
the classical parabolic-elliptic Keller-Segel system obtained from (\ref{1}) in the limit case $\Omega=\R^2$.
Indeed, choosing $a\equiv 4$ and $b\equiv const.$ in (\ref{foo1}) one would rediscover a well-known
family of explicit radial steady states for the corrseponding version of (\ref{1}) (\cite{lopezg_nagai_yamada}).
\mysection{Local existence}
The following basic result on local existence of solutions to (\ref{0})
can be proved by adapting approaches that are well-established in the context of
parabolic-elliptic models for taxis mechanisms involving both cross-diffusion terms and
ODE dynamics (cf.~\cite{taowin_proca}, \cite{meral_stinner_surulescu}, \cite{litcanu_morales-rodrigo}
and \cite{ciewin}, for instance).
Here we note that our assumption that $w_0$ belong to $C^1(\overline{\Omega})$ enables us to
use standard elliptic Schauder theory to gain appropriate knowledge on
the spatial regularity of $v$. Indeed, expressing $w$ via the formula
\be{w_explicit}
    w(x, t) =w_0(x) e^{-\frac{\delta}{\tau} t} + \frac{1}{\tau} \int_0^t e^{-\frac{\delta}{\tau} (t-s)} u(x, s) ds,
        \qquad x\in \Omega, \ t>0,
\ee
we see that $v(\cdot,t)$ actually solves the Poisson equation with a temporally nonlocal inhomogeneity
which thanks to the inclusion $w_0\in C^1(\overline{\Omega})$
will be H\"older continuous in $\overline{\Omega}$ provided that $u(\cdot,t)$ is sufficiently regular, where the
latter can be guaranteed by standard arguments involving appropriate smoothing properties of the Neumann heat
semigroup in $\Omega$.
\begin{lem}\label{lem_loc}
  Let $\delta\ge 0$, and suppose that
  $u_0\in C^0(\bar{\Omega})$ and $w_0\in C^1(\bar{\Omega})$ are nonnegative.
  Then there exist $\tm \in (0,\infty]$ and uniquely determined nonnegative functions
  \bas
    & & u\in C^0(\overline{\Omega}\times [0,\tm)) \cap C^{2,1}(\overline{\Omega}\times (0,\tm)), \\
    & & v\in C^{2,0}(\overline{\Omega}\times [0,\tm)), \\
    & & w \in C^{0,1}(\overline{\Omega}\times [0,\tm)),
  \eas
  which solve (\ref{0}) classically in $\Omega \times (0,\tm)$ and which are such that
  \be{extend}
    \mbox{if $\tm<\infty$, then}~
    \|u(\cdot, t)\|_{L^\infty(\Omega)}\rightarrow
    \infty \qquad \mbox{as $ t\nearrow \tm $}.
  \ee
\end{lem}
The following identities describing the evolution of the total masses of the first and third components in (\ref{0})
can easily be checked.
\begin{lem}\label{basic}
  Let $\delta\ge 0$. Then the solution $(u,v,w)$ of (\ref{0}) satisfies
  \be{mass}
    \io u(\cdot,t)=m:=\io u_0 \qquad \mbox{for all } t \in (0,\tm),
  \ee
  and for all $t\in (0,\tm)$ we have
  \bea{mass_w}
    \io w(\cdot,t) &=& e^{-\frac{\delta}{\tau} t} \io w_0 + \frac{m}{\tau} \int_0^t e^{-\frac{\delta}{\tau}(t-s)} ds
        \nn\\[1mm]
    &=& \left\{ \begin{array}{ll}
        e^{-\frac{\delta}{\tau} t} \io w_0 + \frac{m}{\delta}(1-e^{-\frac{\delta}{\tau} t})
        \qquad & \mbox{if } \delta>0, \\[1mm]
        e^{-\frac{\delta}{\tau} t} \io w_0 + \frac{m}{\tau} \cdot t & \mbox{if } \delta=0,
    \end{array} \right.
  \eea
\end{lem}
\proof
  Integrating the first equation in (\ref{0}) with respect
  to $x\in\Omega$, we see that $\frac{d}{dt}\io u\equiv 0$, which immediately yields (\ref{mass}).
  Using this, we only need to integrate (\ref{w_explicit}) in space to obtain (\ref{mass_w}).
\qed
Based on (\ref{mass_w}) we can now explicitly rewrite the degradation term $\mu(t)$ in the second equation in (\ref{0}).
\begin{cor}\label{cor_mu}
  Let $\delta\ge 0$. Then the function $\mu$ defined in (\ref{mu_def}) is given by
  \be{mu}
    \mu(t)=\frac{1}{\pi} e^{-\frac{\delta}{\tau} t} \io w_0 + \frac{m}{\pi\tau} \int_0^t e^{-\frac{\delta}{\tau}(t-s)} ds
    \qquad \mbox{for all } t\in (0,\tm),
  \ee
  where $m:=\io u_0$.
\end{cor}
\mysection{Global existence}
The following basic statement on the time evolution of the functional $\frac{1}{p} \io u^p + \frac{\tau}{p+1} \io w^{p+1}$
will be the starting point for our derivation of bounds for $u$, and also for $w$, in spaces of the form
$L^\infty ((0, \tm); L^p(\Omega))$ with $p>1$. Besides in Lemma \ref{lem22}, it will be referred to in Lemma \ref{lem25}
below.
\begin{lem}\label{lem21}
  Let $\delta\ge 0$. Then for all $p>1$, the solution of (\ref{0}) satisfies
  \be{21.1}
    \frac{d}{dt} \bigg\{ \frac{1}{p} \io u^p + \frac{\tau}{p+1} \io w^{p+1} \bigg\}
    + \frac{4(p-1)}{p^2} \io |\nabla u^\frac{p}{2}|^2 + \delta \io w^{p+1}
    \le \frac{p-1}{p} \io u^p w + \io uw^p
  \ee
  for all $t\in (0,\tm)$.
\end{lem}
\proof
  We multiply the first equation in (\ref{0}) by $u^{p-1}$ and
  integrate by parts using the identity $\Delta v=\mu(t)-w$ to find
  that
  \bea{21.01}
    \frac{1}{p}\frac{d}{dt}\io u^p +\frac{4(p-1)}{p^2}\io |\nabla u^{\frac{p}{2}}|^2
    &=& (p-1) \io u^{p-1}\nabla u\cdot \nabla v\nn\\[1mm]
    &=&-\frac{p-1}{p} \io u^p\Delta v\nn\\[1mm]
    &=&-\frac{p-1}{p} \io u^p \Big(\mu(t)-w\Big)\nn\\[1mm]
    &\le& \frac{p-1}{p} \io u^p w \qquad \mbox{for all $t\in (0,\tm)$},
  \eea
  because $\mu(t)\ge 0$ by Corollary \ref{cor_mu}. On the other
  hand, multiplying the third equation in (\ref{0}) by $w^p$ and
  integrating with respect to $x\in \Omega$ we see that
  \bas
    \frac{\tau}{p+1}\frac{d}{dt}\io w^{p+1} +\delta \io w^{p+1} =\io uw^p
    \qquad \mbox{for all $t\in (0,\tm)$}.
  \eas
  Adding this to (\ref{21.01}) proves (\ref{21.1}).
\qed
Further estimating the terms on the right of (\ref{21.1})
shows that the functional in question actually satisfies the following autonomous differential inequality.
\begin{lem}\label{lem22}
  Let $\delta\ge 0$. Then
  for any $p>1$ there exists $C(p)>0$ such that the solution of (\ref{0}) satisfies
  \be{22.1}
    \frac{d}{dt} \bigg\{ \frac{1}{p} \io u^p + \frac{\tau}{p+1} \io w^{p+1} \bigg\}
    \le C(p) \cdot \bigg\{ \frac{1}{p} \io u^p + \frac{\tau}{p+1} \io w^{p+1} \bigg\}
    \qquad \mbox{for all } t\in (0,\tm).
  \ee
\end{lem}
\proof
  Let us first invoke the Gagliardo-Nirenberg inequality to fix $c_1>0$ such that
  \bea{22.11}
    \io \varphi^{p+1}
    &=& \|\varphi^{\frac{p}{2}}\|_{L^{\frac{2(p+1)}{p}}(\Omega)}^{\frac{2(p+1)}{p}}\nn\\[1mm]
    &\le& c_1 \|\nabla \varphi^{\frac{p}{2}}\|_{L^2(\Omega)}^2
        \cdot \|\varphi^{\frac{p}{2}}\|_{\frac{2}{p}(\Omega)}^{\frac{2}{p}}
        + c_1 \|\varphi^{\frac{p}{2}}\|_{\frac{2}{p}(\Omega)}^{\frac{2(p+1)}{p}}\nn\\[1mm]
    &=& c_1 \|\nabla \varphi^{\frac{p}{2}}\|_{L^2(\Omega)}^2 \cdot \|\varphi\|_{L^1(\Omega)}
        + c_1\|\varphi\|_{L^1(\Omega)}^{p+1} \qquad
    \mbox{for all nonnegative $\varphi \in W^{1,2}(\Omega)$.}
  \eea
  We now let $\eps:=\frac{2(p-1)}{mc_1 p^2}$
  and use the Young inequality to estimate the two terms on the right of (\ref{21.1}) according to
  \be{22.01}
    \frac{p-1}{p} \io u^p w + \io uw^p \le 2\eps \io u^{p+1} +(\eps^{-p}
    +\eps^{-\frac{1}{p}}) \io w^{p+1} \qquad \mbox{for all $t\in (0, \tm)$.}
  \ee
  Here since $\|u\|_{L^1(\Omega)}=\io u=m$ for all $t\in (0, \tm)$ due to Lemma \ref{basic},
  by the H\"{o}lder inequality and (\ref{22.11}) we obtain
  \bas
    2\eps \io u^{p+1} &\le& 2\eps c_1 m \io |\nabla u^{\frac{p}{2}}|^2
    +2\eps c_1 m \bigg(\io u\bigg)^p\nn \\[1mm]
    &\le& 2\eps c_1 m \io |\nabla u^{\frac{p}{2}}|^2 +2\eps c_1 m\cdot |\Omega|^\frac{p-1}{p} \io u^p \\[1mm]
    &=& \frac{4(p-1)}{p^2} \io |\nabla u^{\frac{p}{2}}|^2 + \frac{4(p-1)}{p^2} \cdot |\Omega|^\frac{p-1}{p} \io u^p
    \qquad \mbox{for all } t\in (0,\tm).
  \eas
  Inserting this into (\ref{22.01}) and recalling (\ref{21.1}) proves (\ref{22.1}).
\qed
We are now in the position to assert our global existence result for (\ref{0}).\abs
\proofc of Proposition \ref{prop100}.\quad
  For any given $T\in (0, \tm)$, the ODI (\ref{22.1}) yields
  \bas
    \frac{1}{p} \io u^p +\frac{\tau}{p+1}\io w^{p+1}\le c_1(p, T)
    \qquad\mbox{for all $t\in (0, T)$}
  \eas
  with $c_1(p, T) :=\Big(\frac{1}{p} \io u_0^p +\frac{\tau}{p+1}\io
  w_0^{p+1}\Big) \cdot e^{C(p)\cdot T}$, where $C(p)>0$ is as defined by
  Lemma \ref{lem22}. Since $\tau>0$, this immediately yields
  \be{100.02}
    \io u^p \le p c_2(p, T) \qquad\mbox{for all $t\in (0, T)$}
  \ee
  and
  \bas
    \io w^{p+1} \le c_2(p, T) \qquad\mbox{for all $t\in (0, T)$}
  \eas
  with $c_2(p, T) :=\max\{p,\frac{p+1}{\tau}\} \cdot c_1(p, T)$. From the latter and standard
  elliptic regularity theory we obtain a bound for $v$ in all spaces
  $L^\infty((0, T); W^{2, p}(\Omega))$ for any $p\in (1, \infty)$,
  whence in particular there exists $c_3(p,T)>0$ such that
  \bas
    \|\nabla v(\cdot, t)\|_{L^\infty(\Omega)}\le c_3(p, T) \qquad \mbox{for all $t\in (0, T)$.}
  \eas
  Along with (\ref{100.02}), this ensures
  that Lemma 4.1 in \cite{taowin2} becomes applicable so as to assert
  via a Moser-type iteration that
  \bas
    \|u(\cdot, t)\|_{L^\infty(\Omega)} \le c_4(p, T) \qquad\mbox{for all $t\in (0, T)$}
  \eas
  holds for some $c_4(p, T)>0$. Finally, Proposition \ref{prop100} is an evident consequence of this and the
  extensibility criterion in Lemma \ref{lem_loc}.
\qed
\mysection{Radial solutions. A comparison principle}
Throughout the sequel we shall assume that the initial data $u_0$ and $w_0$, and hence clearly also all components of
the solution $(u,v,w)$, are radially symmetric with respect to the spatial origin, and unless stated otherwise we fix
\be{m}
    m:=\io u_0.
\ee
Then without danger of confusion we may and will switch to the usual radial notation and write $u=u(r,t)$ for
$r=|x| \in [0,1]$ whenever this appears convenient.
\begin{lem}\label{lem31}
  Suppose that $\delta\ge 0$, and given a radial solution $(u,v,w)$ of (\ref{0}), let
  \be{31.1}
    U(\xi,t):=\int_0^{\sqrt{\xi}} ru(r,t)dr,
    \qquad \xi \in [0,1], \ t\ge 0.
  \ee
  Then
  \be{31.01}
    U(0,t)=0
    \quad \mbox{as well as} \quad
    U(1,t)=\frac{m}{2\pi}
    \qquad \mbox{for all } t\ge 0
  \ee
  and
  \be{31.02}
    U_\xi(\xi,t) \ge 0
    \qquad \mbox{for all $\xi\in (0,1)$ and } t>0.
  \ee
  Moreover,
  \be{31.11}
    \parab U(\xi,t) =0
    \qquad \mbox{for all $\xi\in (0,1)$ and } t>0,
  \ee
  where the operator $\parab$ is defined by setting
  \bea{parab}
    \hspace*{-3mm}
    \parab \tU(\xi,t) :=
    \tU_t - 4\xi \tU_{\xi\xi}
    - \frac{2}{\tau} \Bigg\{ \int_0^t e^{-\frac{\delta}{\tau}(t-s)}
        \Big( \tU(\xi,s)-\frac{m}{2\pi} \xi\Big) \, ds \Bigg\} \cdot \tU_\xi
    - 2\Big(W_0(\xi)-\kz \xi\Big) \cdot e^{-\frac{\delta}{\tau} t} \tU_\xi
  \eea
  for $\xi\in (0,1), t>0$ and $\tU\in C^1((0,1)\times (0,\infty)) \cap C^0((0,\infty);W^{2,\infty}((0,1)))$, with
  \be{W0}
    W_0(\xi):=\int_0^{\sqrt{\xi}} rw_0(r)dr, \quad \xi\in [0,1],
    \qquad \mbox{and} \qquad
    \kz:=W_0(1)=\int_0^1 rw_0(r)dr.
  \ee
\end{lem}
\proof
  The boundary properties in (\ref{31.01}) are immediate from (\ref{31.1}) and (\ref{mass}), whereas the monotonicity
  statement in (\ref{31.02}) is equivalent to the nonnegativity of $u$.
  Moreover, upon differentiation in (\ref{0}) we see that for $\xi\in (0,1)$ and $t>0$,
  \bas
    U_t(\xi,t) &=& \int_0^{\sqrt{\xi}} r\cdot \Big\{ \frac{1}{r} (ru_r)_r - \frac{1}{r} (ruv_r)_r \Big\} \, dr \\
    &=& \sqrt{\xi} u_r(\sqrt{\xi},t) - \sqrt{\xi} u(\sqrt{\xi},t) v_r(\sqrt{\xi},t),
  \eas
  where by (\ref{31.1}) we have
  \bas
    u(\sqrt{\xi},t)=2U_\xi(\xi,t)
    \qquad \mbox{and} \qquad
    u_r(\sqrt{\xi},t)=4\sqrt{\xi} U_{\xi\xi}(\xi,t).
  \eas
  Since the second equation in (\ref{0}) implies that
  \bas
    rv_r(r,t)=-\int_0^r \rho w(\rho,t)d\rho + \frac{\mu(t)r^2}{2}
    \qquad \mbox{for all $r\in (0,1)$ and } t>0,
  \eas
  we thus obtain
  \be{31.2}
    U_t=4\xi U_{\xi\xi} + 2U_\xi W - \mu(t) \xi U_\xi
    \qquad \mbox{for all $\xi\in (0,1)$ and } t>0,
  \ee
  with $W(\xi,t):=\int_0^{\sqrt{\xi}} rw(r,t)dr$, $\xi\in [0,1], t\ge 0$.
  Now by (\ref{W0}) and (\ref{w_explicit}),
  \bas
    W(\xi,t)=W_0(\xi) e^{-\frac{\delta}{\tau} t} + \frac{1}{\tau} \int_0^t e^{-\frac{\delta}{\tau}(t-s)} U(\xi,s)ds
    \qquad \mbox{for all $\xi\in (0,1)$ and } t>0,
  \eas
  whereas
  \bas
    \mu(t)=2\kz e^{-\frac{\delta}{\tau} t} + \frac{m}{\pi\tau} \int_0^t e^{-\frac{\delta}{\tau} (t-s)} ds
    \qquad \mbox{for all } t>0
  \eas
  according to (\ref{mu}) and (\ref{W0}).
  Therefore,
  \bas
    2U_\xi W - \mu(t)\xi U_\xi
    &=& \frac{2}{\tau} \bigg\{ \int_0^t e^{-\frac{\delta}{\tau}(t-s)} U(\xi,s)ds \bigg\} \cdot U_\xi(\xi,t)
    + 2W_0(\xi) e^{-\frac{\delta}{\tau} t} U_\xi(\xi,t) \\
    & & - 2\kz \xi \cdot e^{-\frac{\delta}{\tau} t} U_\xi(\xi,t)
    - \bigg\{ \int_0^t e^{-\frac{\delta}{\tau}(t-s)} \cdot \frac{m}{\pi\tau} \xi ds \bigg\} \cdot U_\xi(\xi,t)
  \eas
  for all $\xi\in (0,1)$ and $t>0$, which along with (\ref{31.2}) proves (\ref{31.11}).
\qed
Fortunately, the parabolic operator $\parab$ introduced above falls into a class of operators allowing for
a comparison principle.
To see this, for functions $A,B$ and $D$ to be specified below, let us consider
\be{qarab}
    \qarab \tU(\xi,t)
    := \tU_t(\xi,t)-A(\xi,t) \tU_{\xi\xi}(\xi,t)
    - \bigg\{B(\xi,t)+ \int_0^t D(\xi,t,s) \tU(\xi,s)ds \bigg\} \cdot \tU_\xi(\xi,t),
    \quad \xi\in (0,1), \ t\in (t_0,T),
\ee
for $0\le t_0<T$ and sufficiently regular $\tU:(0,1)\times (0,T) \to \R$.
Then assuming, besides parabolicity, that the memory term has a favorable sign, we can indeed derive the following
comparison principle for spatially nondecreasing functions.
\begin{lem}\label{lem12}
  Let $t_0\ge 0$ and $T>t_0$, and suppose that $A\in C^0((0,1) \times (t_0,T))$, $B\in C^0((0,1)\times (t_0,T))$ and
  $D\in C^0([0,1] \times [0,T] \times [0,T])$ satisfy
  \be{12.1}
    A \ge 0 \quad \mbox{in } (0,1) \times (t_0,T)
    \qquad \mbox{and} \qquad
    D\ge 0 \quad \mbox{in } [0,1] \times [0,T] \times [0,T].
  \ee
  Moreover, assume that $\UU$ and $\OU$ are nonnegative functions belonging to
  \be{12.11}
    C^0([0,1]\times [0,T]) \cap C^1((0,1)\times (t_0,T)) \cap C^0((t_0,T);W^{2,\infty}((0,1))),
  \ee
  which are such that
  \be{12.2}
    0 \le \UU_\xi(\xi,t) \le M
    \qquad \mbox{for all $\xi\in (0,1)$ and } t\in (t_0,T)
  \ee
  with some $M>0$, and such that with $\qarab$ as defined in (\ref{qarab}) we have
  \be{12.3}
    \qarab \UU(\cdot,t) \le \qarab \OU(\cdot,t)
    \qquad \mbox{a.e.~in $(0,1)$ for all } t\in (t_0,T).
  \ee
  Then if
  \be{12.4}
    \UU(\xi,t)\le \OU(\xi,t)
    \qquad \mbox{for all $\xi\in [0,1]$ and } t\in [0,t_0]
  \ee
  as well as
  \be{12.5}
    \UU(0,t)\le \OU(0,t)
    \quad \mbox{for all } t\in [t_0,T]
    \qquad \mbox{and} \qquad
    \UU(1,t)\le \OU(1,t)
    \quad \mbox{for all } t\in [t_0,T],
  \ee
  we have the global ordering property
  \be{12.6}
    \UU(\xi,t) \le \OU(\xi,t)
        \quad \mbox{for all $\xi\in [0,1]$ and } t\in [0,T].
  \ee
\end{lem}
\proof
  We let $c_1:=\|D\|_{L^\infty((0,1)\times (0,T)\times (0,T))}$ and $\alpha:=\sqrt{c_1 M}$ with $M$ as in (\ref{12.2}),
  and for arbitrary $\eps_0>0$ we let
  \be{12.61}
    \eps(t):=\eps_0 e^{\alpha t}, \qquad t\ge 0,
  \ee
  and
  \bas
    d(\xi,t):=\UU(\xi,t)-\OU(\xi,t)-\eps(t)
    \qquad \mbox{for $\xi\in [0,1]$ and } t\in [0,T].
  \eas
  Then $d$ is continuous in $[0,1]\times [0,T]$ with
  \bas
    d(\xi,t) \le -\eps_0 e^{\alpha t} < 0
    \qquad \mbox{for all $\xi\in [0,1]$ and } t\in [0,t_0]
  \eas
  by (\ref{12.4}) and
  \bas
    d(\xi,t) \le -\eps_0 e^{\alpha t} < 0
    \qquad \mbox{for $\xi\in \{0,1\}$ and } t\in [t_0,T]
  \eas
  according to (\ref{12.5}).
  Thus,
  \bas
    t_\star:=\sup \Big\{ t\in (0,T) \ \Big| \ d<0 \mbox{ in } [0,1] \times [0,t] \Big\}
  \eas
  is well-defined and satisfies $t_\star\in (t_0,T]$, and if we had $t_\star<T$, then there would exist
  $\xi_\star\in (0,1)$ such that
  \be{12.7}
    d(\xi_\star,t_\star)=\max_{\xi\in [0,1]} d(\xi,t_\star)=0,
  \ee
  whence evidently
  \be{12.8}
    d_t(\xi_\star,t_\star) \ge 0
    \qquad \mbox{and} \qquad
    d_\xi(\xi_\star,t_\star)=0,
  \ee
  because $d\in C^1((0,1)\times (t_0,T))$ by (\ref{12.11}).
  Now by (\ref{12.3}) we know that there exists a null set $N\subset (0,1)$ such that $d_{\xi\xi}(\xi,t_\star)$ exists
  for all $\xi \in (0,1)\setminus N$ and
  \bea{12.9}
    d_t(\xi,t_\star)
    &\le& A(\xi,t_\star) d_{\xi\xi}(\xi,t_\star)
    + \bigg\{ B(\xi,t_\star) + \int_0^{t_\star} d(\xi,t_\star,s) \OU(\xi,s)ds \bigg\} \cdot d_\xi(\xi,t_\star) \nn\\[1mm]
    & & + \UU_\xi(\xi,t_\star) \cdot \int_0^{t_\star} D(\xi,t_\star,s) d(\xi,s) ds \nn\\[1mm]
    & & + \UU_\xi(\xi,t_\star) \cdot \int_0^{t_\star} D(\xi,t_\star,s) \cdot \eps(s)ds
    - \eps'(t_\star)
    \qquad \mbox{for all } \xi\in (0,1)\setminus N.
  \eea
  In order to make appropriate use of (\ref{12.8}) and the maximality property in (\ref{12.7}), we observe that (\ref{12.7})
  necessarily implies that there exists $(\xi_j)_{j\in\N} \subset (0,1)\setminus N$ such that
  $\xi_j\to \xi_\star$ as $j\to\infty$ and
  \bas
    d_{\xi\xi}(\xi_j,t_\star) \le 0
    \qquad \mbox{for all }j\in\N,
  \eas
  for otherwise we would have $\rm{essliminf}_{\xi\to\xi_\star} d_{\xi\xi}(\xi,t_\star)>0$, contradicting (\ref{12.7}).
  Choosing $\xi=\xi_j$ in (\ref{12.9}) and using that (\ref{12.8}) and (\ref{12.1}) entail that
  \bas
    \limsup_{j\to\infty} d_t(\xi_j,t_\star)\ge 0
    \qquad \mbox{and} \qquad
    d_\xi(\xi_j,t_\star)\to 0
    \qquad \mbox{as } j\to\infty,
  \eas
  we obtain on taking $j\to\infty$ that
  \bea{12.10}
    0 &\le& \UU_\xi(\xi_\star,t_\star) \cdot \int_0^{t_\star} D(\xi_\star,t_\star,s) d(\xi_\star,s)ds
    + \UU_\xi(\xi_\star,t_\star) \cdot \int_0^{t_\star} D(\xi_\star,t_\star,s) \cdot \eps(s)ds
    -\eps'(t_\star).
  \eea
  Here since $d(\xi_\star,s)<0$ for all $s\in [0,t_\star)$ by definition of $t_\star$, we have
  \bas
    \UU_\xi(\xi_\star,t_\star) \cdot \int_0^{t_\star} D(\xi_\star,t_\star,s) d(\xi_\star,s)ds\le 0,
  \eas
  because $D\ge 0$ and $\UU_\xi \ge 0$ by (\ref{12.1}) and (\ref{12.2}).
  Furthermore, (\ref{12.2}) and our choice of $c_1$ ensure that
  \bas
    \UU_\xi(\xi_\star,t_\star) \cdot \int_0^{t_\star} D(\xi_\star,t_\star,s) \cdot \eps(s)ds
    \le Mc_1 \int_0^{t_\star} \eps(s)ds,
  \eas
  so that recalling (\ref{12.61}), from (\ref{12.10}) we obtain
  \bas
    0 &\le& Mc_1 \int_0^{t_\star} \eps_0 e^{\alpha s} ds - \alpha \eps_0 e^{\alpha t_\star} \\
    &=& \frac{Mc_1 \eps_0}{\alpha} (e^{\alpha t_\star}-1) - \alpha \eps_0 e^{\alpha t_\star} \\
    &<& \Big(\frac{Mc_1}{\alpha^2}-1\Big) \cdot \alpha \eps_0 e^{\alpha t_\star}.
  \eas
  Since $\frac{Mc_1}{\alpha^2}=1$ according to our definition of $\alpha$, this absurd conclusion shows that actually
  $t_\star=T$ and hence $\UU \le \OU + \eps_0 e^{\alpha t}$ throughout $[0,1] \times [0,T]$. On taking
  $\eps_0\searrow 0$ we finally arrive at the desired inequality.
\qed
\mysection{Boundedness for $\io u_0<8\pi\delta$. Proof of Theorem \ref{theo200}}
In this section we shall make sure that small-mass solutions remain bounded in the sense of Theorem \ref{theo200}.
\subsection{A pointwise upper bound for $U$}
As a preliminary, let us prove the following elementary lemma.
\begin{lem}\label{lem9}
  Let $m>0$ and $\eps>0$, and suppose that $\varphi\in W^{1,\infty}((0,1))$ is such that $\varphi(0)=0$ and
  $\varphi(\xi) \le \frac{m}{2\pi}$ for all $\xi\in (0,1)$. Then there exists $b\in (0,1)$ such that
  \be{9.1}
    \varphi(\xi) \le \frac{m}{2\pi} \cdot \frac{(b+1+\eps)\xi}{b+\xi}
    \qquad \mbox{for all } \xi\in (0,1).
  \ee
\end{lem}
\proof
  Since $\varphi(0)=0$ and $\varphi_\xi\in L^\infty(\Omega)$, we can find $c_1>0$ such that
  $\varphi(\xi) \le c_1\xi$ for all $\xi\in (0,1)$, where we may assume that $c_1>\frac{m}{2\pi}$.
  Therefore, our assumption warrants that
  \be{9.2}
    \varphi(\xi) \le \hat\varphi(\xi):=\min \Big\{ c_1\xi , \frac{m}{2\pi} \Big\}
    \qquad \mbox{for all } \xi\in (0,1).
  \ee
  Now writing $\varphi_b(\xi):=\frac{m}{2\pi} \cdot \frac{(b+1+\eps)\xi}{b+\xi}$ for $\xi\in [0,1]$ and
  $b\in (0,1)$, we see that
  the quotient $\frac{\hat\varphi}{\varphi_b}$ admits a continuous extension $Q_b$ to all of $[0,1]$ such that
  \bas
    Q_b(\xi)=\left\{ \begin{array}{ll}
    \frac{2\pi}{m} \cdot \frac{c_1(b+\xi)}{b+1+\eps} \qquad & \mbox{if } \xi \in [0,\xi_1], \\[1mm]
    \frac{b+\xi}{(b+1+\eps)\xi} & \mbox{if } \xi\in (\xi_1,1],
    \end{array} \right.
  \eas
  where $\xi_1:=\frac{m}{2\pi c_1} \in (0,1)$ thanks to our choice of $c_1$.
  Since
  \bas
    \frac{\partial}{\partial b} \frac{b+\xi}{b+1+\eps}=\frac{1+\eps-\xi}{(b+1+\eps)^2}>0
    \qquad \mbox{for all } \xi\in [0,1],
  \eas
  it follows that as $b\searrow 0$ we have
  \bas
    Q_b(\xi)\searrow Q(\xi):=\left\{ \begin{array}{ll}
    \frac{2\pi}{m} \cdot \frac{c_1\xi}{1+\eps} \qquad & \mbox{if } \xi \in [0,\xi_1], \\[1mm]
    \frac{1}{1+\eps} & \mbox{if } \xi\in (\xi_1,1].
    \end{array} \right.
  \eas
  As $Q$ is continuous in $[0,1]$, Dini's theorem asserts that the convergence $Q_b\to Q$ is actually uniform in $[0,1]$.
  Since $Q(\xi)\le \frac{1}{1+\eps}<1$ for all $\xi\in [0,1]$, we can therefore pick some sufficiently small $b\in (0,1)$
  such that $Q_b(\xi)\le 1$ for all $\xi\in [0,1]$, which in view of (\ref{9.2}) implies (\ref{9.1}).
\qed
By means of a comparison argument, we can now prove that under the assumption $\io u_0<8\pi\delta$, it is possible
to control the mass concentrating in small balls around the origin uniformly with respect to $t\in (0,\infty)$
in the following sense.
\begin{lem}\label{lem10}
  Let $\delta\ge 0$, and assume that $u_0$ has the property that
  \be{10.1}
    m\equiv \io u_0<8\pi\delta.
  \ee
  Then there exists $C>0$ such that the function $U$ defined in (\ref{31.1}) satisfies
  \be{10.2}
    U(\xi,t)\le C\xi
    \qquad \mbox{for all $\xi\in (0,1)$ and } t>0.
  \ee
\end{lem}
\proof
  Since $m<8\pi\delta$, we can find $\eps>0$ such that
  \bas
    8>\frac{m}{\pi\delta} (1+\eps),
  \eas
  and thereupon choose $t_0>0$ large fulfilling
  \be{10.3}
    8\ge \frac{m}{\pi\delta}(1+\eps) + 4c_1 e^{-\frac{\delta}{\tau} t_0},
  \ee
  where
  \be{10.33}
    c_1:=\frac{1}{2} \|w_0\|_{L^\infty(\Omega)}.
  \ee
  With these values of $\eps$ and $t_0$ fixed, using that for all $t\in [0,t_0]$ we have $U(0,t)=0$ and
  $U(\xi,t)\le U(1,t) =\frac{m}{2\pi}$ for all $\xi\in [0,1]$ by (\ref{31.01}) and (\ref{31.02}),
  we can apply Lemma \ref{lem9} to find $b\in (0,1)$ satisfying
  \be{10.4}
    \max_{t\in [0,t_0]} U(\xi,t) \le \frac{m}{2\pi} \cdot \frac{(b+1+\eps)\xi}{b+\xi}
    \qquad \mbox{for all } \xi\in [0,1].
  \ee
  This means that if we let
  \bas
    \OU(\xi,t):=\frac{a\xi}{b+\xi}
    \qquad \mbox{for $\xi\in [0,1]$ and } t\ge 0,
  \eas
  with
  \bas
    a:=\frac{m}{2\pi} \cdot (b+1+\eps),
  \eas
  then $U(\xi,t)\le \OU(\xi,t)$ for all $\xi\in [0,1]$ and $t\in [0,t_0]$. Moreover, clearly
  $0=U(0,t)\le \OU(0,t)=0$ and $U(1,t)=\frac{m}{2\pi}<\frac{m}{2\pi} \cdot \frac{b+1+\eps}{b+1}=\OU(1,t)$ for all
  $t\ge t_0$.
  Computing
  \bas
    \OU_t=0,\quad  \OU_\xi=\frac{ab}{(b+\xi)^2}
    \quad \mbox{and} \quad
    \OU_{\xi\xi}=- \frac{2ab}{(b+\xi)^3}
    \qquad \mbox{for $\xi\in (0,1)$ and } t>t_0,
  \eas
  we see that with $\parab$ as defined in (\ref{parab}) we have
  \bea{10.5}
    \parab \OU(\xi,t)
    &=& \frac{8ab\xi}{(b+\xi)^3}
    - \bigg\{ \frac{2}{\tau}\int_0^t e^{-\frac{\delta}{\tau}(t-s)} \cdot
        \Big(\frac{a\xi}{b+\xi} - \frac{m}{2\pi} \xi\Big) ds \bigg\}
    \cdot \frac{ab}{(b+\xi)^2} \nn\\
    & & - 2\Big(W_0(\xi)-\kz \xi\Big) \cdot e^{-\frac{\delta}{\tau} t} \cdot \frac{ab}{(b+\xi)^2} \nn\\[2mm]
    &=& \frac{ab\xi}{(b+\xi)^2} \cdot \bigg\{
    \frac{8}{b+\xi} - \frac{2}{\delta} (1-e^{-\frac{\delta}{\tau} t}) \cdot \Big(\frac{a}{b+\xi} - \frac{m}{2\pi}\Big)
    - 2\Big(\frac{W_0(\xi)}{\xi} - \kz\Big) \cdot e^{-\frac{\delta}{\tau} t} \bigg\}
  \eea
  for all $\xi\in (0,1)$ and $t>t_0$.
  Here we use the definition of $a$ and the nonnegativity of $e^{-\frac{\delta}{\tau} t}$ to estimate
  \bas
    \frac{2}{\delta}(1-e^{-\frac{\delta}{\tau} t}) \cdot \Big(\frac{a}{b+\xi} - \frac{m}{2\pi}\Big)
    &=& \frac{2}{\delta}(1-e^{-\frac{\delta}{\tau} t})\cdot \frac{m}{2\pi} \cdot \Big(\frac{b+1+\eps}{b+\xi}-1\Big) \\
    &=& \frac{2}{\delta}(1-e^{-\frac{\delta}{\tau} t})\cdot \frac{m}{2\pi} \cdot \frac{1+\eps-\xi}{b+\xi} \\
    &\le& \frac{2}{\delta} \cdot \frac{m}{2\pi} \cdot \frac{1+\eps}{b+\xi}
    \qquad \mbox{for all $\xi\in (0,1)$ and } t>t_0.
  \eas
  Since by (\ref{W0}) and (\ref{10.33}) we have
  \bas
    W_0(\xi) \le \|w_0\|_{L^\infty(\Omega)} \cdot \frac{(\sqrt{\xi})^2}{2} = c_1\xi
    \qquad \mbox{for all } \xi\in (0,1),
  \eas
  we moreover see that
  \bas
    2\Big(\frac{W_0(\xi)}{\xi}-\kz\Big) \cdot e^{-\frac{\delta}{\tau} t}
    &\le& 2c_1 e^{-\frac{\delta}{\tau} t} \le 2c_1 e^{-\frac{\delta}{\tau} t_0}
    \qquad \mbox{for all $\xi\in (0,1)$ and } t>t_0.
  \eas
  According to (\ref{10.3}), the identity (\ref{10.5}) thus shows that
  \bas
    \parab \OU(\xi,t)
    &\ge& \frac{ab\xi}{(b+\xi)^2} \cdot \bigg\{
    \frac{8}{b+\xi} - \frac{2}{\delta} \cdot \frac{m}{2\pi} \cdot \frac{1+\eps}{b+\xi}
    - 2c_1 e^{-\frac{\delta}{\tau} t_0} \bigg\} \\
    &=& \frac{ab\xi}{(b+\xi)^3} \cdot \bigg\{
    8-\frac{m}{\pi\delta} (1+\eps) - 2c_1 (b+\xi) e^{-\frac{\delta}{\tau} t_0} \bigg\} \\
    &\ge& \frac{ab\xi}{(b+\xi)^3} \cdot \bigg\{
    8-\frac{m}{\pi\delta} (1+\eps) - 4c_1 e^{-\frac{\delta}{\tau} t_0} \bigg\} \\[2mm]
    &\ge& 0
    \qquad \mbox{for all $\xi\in (0,1)$ and } t>t_0,
  \eas
  where we have used that $b+\xi\le b+1\le 2$, because $b<1$.
  By comparison on the basis of Lemma \ref{lem12}, we thereby conclude that $\OU \ge U$ in $(0,1)\times (0,\infty)$,
  which in particular shows that
  \bas
    U(\xi,t) \le \frac{m(b+1+\eps)}{2\pi b} \cdot \xi
    \qquad \mbox{for all $\xi\in (0,1)$ and } t>0,
  \eas
  and thereby completes the proof.
\qed
\subsection{Boundedness away from the origin}
In the case $\delta>0$
when the third equation in (\ref{0}) contains an absorption term, radial solutions can become unbounded in their
first component $u$ only near the spatial origin.
This is contained in the following lemma, the outcome of which will be an essential ingredient to our $\eps$-regularity
result in Section \ref{sect5.3}, and hence in establishing Theorem \ref{theo200}.
\begin{lem}\label{lem24}
  Let $\delta>0$. Then for all $r_0\in (0,1)$ there exists $C(r_0)>0$ such that the solution of (\ref{0}) satisfies
  \be{24.1}
    u(x,t) \le C(r_0)
    \qquad \mbox{for all $x\in \Omega\setminus B_{r_0}$ and } t>0.
  \ee
\end{lem}
\proof
  We evidently only need to consider the case when $u_0\not\equiv
  0$, in which we proceed in six setps.
  \underline{Step 1.} \quad
  We first claim that there exists $c_1>0$ such that
  \be{24.2}
    |v_r(r,t)| \le \frac{c_1}{r}
    \qquad \mbox{for all $r\in (0,1)$ and } t>0.
  \ee
  To verify this, we write the second equation in (\ref{0}) in the form
  \bas
    \frac{1}{r}(rv_r)_r=\mu(t)-w
    \qquad \mbox{for $r\in (0,1)$ and } t>0,
  \eas
  multiply this by $r$ and integrate using $v_r(0,t)=0$ for all $t>0$ to see that
  \bas
    rv_r(r,t)=\frac{\mu(t)r^2}{2} - \int_0^r \rho w(\rho,t)d\rho
    \qquad \mbox{for all $r\in (0,1)$ and } t>0.
  \eas
  Since $w\ge 0$ and $\int_0^1 \rho w(\rho,t)d\rho=\frac{\mu(t)}{2}$ for all $t>0$ by (\ref{mu_def}), from this we obtain
  \bas
    -\frac{\mu(t)}{2} \le rv_r(r,t) \le \frac{\mu(t)r^2}{2}
    \qquad \mbox{for all $r\in (0,1)$ and } t>0.
  \eas
  This implies (\ref{24.2}) if we choose $c_1:=\frac{1}{2}\|\mu\|_{L^\infty((0,\infty))}$ which is finite according to
  (\ref{mu}).\abs
  \underline{Step 2.} \quad
  We next assert that for all $p\in (0,1)$ and each $r_0\in (0,1)$ we can find $c_2(p,r_0)>0$ fulfilling
  \be{24.3}
    \int_t^{t+1} \int_{\Omega\setminus B_{r_0}} |\nabla u^\frac{p}{2}|^2 \le c_2(p,r_0)
    \qquad \mbox{for all } t>0.
  \ee
  To this end, we fix a radially symmetric $\zeta\in C^\infty(\overline{\Omega})$ such that $0\le \zeta\le 1$ in $\Omega$,
  $\zeta\equiv 1$ in $\Omega\setminus B_{r_0}$ and $\zeta\equiv 0$ in $B_{\frac{r_0}{2}}$, and multiply the first equation
  in (\ref{0}) by $\zeta^2 u^{p-1}$ to see upon integrating by parts that
  \bea{24.4}
    \frac{1}{p} \frac{d}{dt} \io \zeta^2 u^p
    &=& (1-p) \io \zeta^2 u^{p-2} |\nabla u|^2
    - 2 \io \zeta u^{p-1} \nabla u \cdot \nabla \zeta \nn\\
    & & - (1-p) \io \zeta^2 u^{p-1} \nabla u \cdot \nabla v
    + 2\io \zeta u^p \nabla v \cdot \nabla \zeta
    \qquad \mbox{for all } t>0,
  \eea
  where we note that $u$, and hence also $u^{p-1}$ and $u^{p-2}$, is smooth and positive in $\overline{\Omega}\times (0,\infty)$
  thanks to our assumption that $u_0\not\equiv 0$ and the strong maximum principle.
  Now by Young's inequality, the H\"older inequality and (\ref{mass}) we have
  \bas
    \bigg| -2\io \zeta u^{p-1} \nabla u \cdot \nabla \zeta \bigg|
    &\le& \frac{1-p}{2} \io \zeta^2 u^{p-2}|\nabla u|^2
    + \frac{1}{1-p} \io u^p |\nabla \zeta|^2 \\
    &\le& \frac{1-p}{2} \io \zeta^2 u^{p-2}|\nabla u|^2
    + \frac{m^p}{1-p} \bigg( \io |\nabla \zeta|^\frac{2}{1-p} \bigg)^{1-p}
    \qquad \mbox{for all } t>0.
  \eas
  By the same token combined with (\ref{24.2}),
  \bas
    \bigg| -(1-p) \io \zeta^2 u^{p-1} \nabla u \cdot \nabla v \bigg|
    &\le& \frac{1-p}{4} \io \zeta^2 u^{p-2} |\nabla u|^2
    + (1-p) \io \zeta^2 u^p |\nabla v|^2 \\
    &\le& \frac{1-p}{4} \io \zeta^2 u^{p-2} |\nabla u|^2
    + (1-p) m^p \bigg(\int_{\Omega\setminus B_\frac{r_0}{2}} |\nabla v|^\frac{2}{1-p} \bigg)^{1-p} \\
    &\le& \frac{1-p}{4} \io \zeta^2 u^{p-2} |\nabla u|^2
    + c_3(p,r_0)
    \qquad \mbox{for all } t>0
  \eas
  with $c_3(p,r_0):=c_1^2 \cdot (1-p) m^p \cdot \Big(2\pi \int_{\frac{r_0}{2}}^1 r^{1-\frac{2}{1-p}} dr \Big)^{1-p}$.
  Similarly, we find $c_4(p,r_0)>0$ fulfilling
  \bas
    \bigg| 2\io \zeta u^p \nabla v \cdot \nabla \zeta \bigg|
    \le 2m^p \bigg( \int_{\Omega \setminus B_{\frac{r_0}{2}}} |\nabla v|^\frac{1}{1-p} \cdot |\nabla \zeta|^\frac{1}{1-p}
        \bigg)^{1-p}
    \le c_4(p,r_0)
    \qquad \mbox{for all } t>0,
  \eas
  whence (\ref{24.4}) altogether yields $c_5(p,r_0)>0$ such that
  \bas
    \frac{1}{p} \frac{d}{dt} \io \zeta^2 u^p
    \ge \frac{1-p}{4} \io \zeta^2 u^{p-2} |\nabla u|^2 - c_5(p,r_0)
    \qquad \mbox{for all } t>0.
  \eas
  After a time integration and another application of the H\"older inequality and (\ref{mass}), we thus obtain
  \bas
    \frac{1-p}{4} \int_t^{t+1} \io \zeta^2 u^{p-2} |\nabla u|^2
    &\le& \frac{1}{p} \io \zeta^2 u^p(\cdot,t+1) + c_5(p,r_0) \\
    &\le& \frac{m^p \cdot \pi^{1-p}}{p} + c_5(p,r_0)
    \qquad \mbox{for all } t>0,
  \eas
  which entails (\ref{24.3}) in view of the fact that $\zeta\equiv 1$ in $\Omega \setminus B_{r_0}$.\abs
  \underline{Step 3.} \quad
  We now make sure that for all $r_0\in (0,1)$ we can find $c_6(r_0)>0$ satisfying
  \be{24.5}
    \int_t^{t+1} \|u(\cdot,s)\|_{L^\infty(\Omega \setminus B_{r_0})} ds \le c_6(r_0)
    \qquad \mbox{for all } t>0.
  \ee
  To this end, let us fix an arbitrary $p\in (0,1)$. Then again by radial symmetry we may combine the one-dimensional
  version of the Gagliardo-Nirenberg inequality with the outcome of Step 2 and (\ref{mass}) to fix positive constants
  $c_7(r_0)$, $c_8(r_0)$ and $c_9(r_0)$ such that
  \bas
    \int_t^{t+1} \|u^\frac{p}{2}(\cdot,s)\|_{L^\infty((r_0,1))}^\frac {2(p+1)}{p} ds
    &\le& c_7(r_0) \int_t^{t+1} \bigg\{
        \|(u^\frac{p}{2})_r(\cdot,s)\|_{L^2((0,1))}^2
        \cdot \|u^\frac{p}{2}(\cdot,s)\|_{L^\frac{2}{p}((r_0,1))}^\frac{2}{p} \\
    & & \hspace*{23mm}
        + \|u^\frac{p}{2}(\cdot,s)\|_{L^\frac{2}{p}((r_0,1))}^\frac{2(p+1)}{p} \bigg\} ds \\
    &\le& c_8(r_0) \int_t^{t+1} \Big\{ \|(u^\frac{p}{2})_r(\cdot,s)\|_{L^2((r_0,1))}^2 + 1 \Big\} ds \\[2mm]
    &\le& c_9(r_0)
    \qquad \mbox{for all } t>0.
  \eas
  Since $\|u^\frac{p}{2}(\cdot,s)\|_{L^\infty((r_0,1))}^\frac{2(p+1)}{p} = \|u(\cdot,s)\|_{L^\infty((r_0,1))}^{p+1}$,
  an application of the H\"older inequality thereupon yields (\ref{24.5}).\abs
  \underline{Step 4.} \quad
  We proceed to show that for any $r_0\in (0,1)$ there exists $c_{10}(r_0)>0$ such that
  \be{24.6}
    \|w(\cdot,t)\|_{L^\infty(\Omega\setminus B_{r_0})} \le c_{10}(r_0)
    \qquad \mbox{for all } t>0.
  \ee
  Indeed, given $t>0$ we write $I_j:=(t-j-1,t-j)\cap (0,\infty)$ for nonnegative integers $j$, and representing
  $w(\cdot,t)$ according to $w(\cdot,t)=e^{-\frac{\delta}{\tau} t} w_0
  + \frac{1}{\tau} \int_0^t e^{-\frac{\delta}{\tau}(t-s)} u(\cdot,s)ds$ we can estimate
  \bas
    \|w(\cdot,t)\|_{L^\infty(\Omega\setminus B_{r_0})}
    &\le& e^{-\frac{\delta}{\tau} t} \|w_0\|_{L^\infty(\Omega)}
    + \frac{1}{\tau} \sum_{j=0}^\infty \int_{I_j} e^{-\frac{\delta}{\tau}(t-s)}
    \|u(\cdot,s)\|_{L^\infty(\Omega\setminus B_{r_0})} ds \\
    &\le& e^{-\frac{\delta}{\tau} t} \|w_0\|_{L^\infty(\Omega)}
    + \frac{1}{\tau} \sum_{j=0}^\infty e^{-\frac{\delta}{\tau} j} \int_{I_j}
    \|u(\cdot,s)\|_{L^\infty(\Omega\setminus B_{r_0})} ds \\
    &\le& e^{-\frac{\delta}{\tau} t} \|w_0\|_{L^\infty(\Omega)}
    + c_6(r_0) \sum_{j=0}^\infty e^{-\frac{\delta j}{\tau}}.
  \eas
  Since the rightmost series converges thanks to our assumption $\delta>0$, this proves (\ref{24.6}).\abs
  \underline{Step 5.} \quad
  Next, we prove that for each $r_0\in (0,1)$ we can fix $c_{11}(r_0)>0$ such that
  \be{24.7}
    \int_t^{t+1} \int_{\Omega\setminus B_{r_0}} |\nabla u|^2 \le c_{11}(r_0)
    \qquad \mbox{for all } t\ge 1.
  \ee
  Since $t\ge 1$, Step 3 allows us to pick $t_0\in (t-1,t)$ such that
  \be{24.8}
    \|u(\cdot,t_0)\|_{L^\infty(\Omega\setminus B_\frac{r_0}{2})} \le c_6\Big(\frac{r_0}{2}\Big).
  \ee
  Then using $\zeta$ as introduced in Step 2, by a straightforward testing procedure we infer that
  \bas
    \frac{1}{2} \frac{d}{dt} \io \zeta u^2(\cdot,s)
    + \io \zeta |\nabla u(\cdot,s)|^2
    &=& - \io u\nabla u \cdot \nabla \zeta
    +\io \zeta u \nabla u\cdot \nabla v
    + \io u^2 \nabla v \cdot \nabla \zeta \\
    &=& \frac{1}{2} \io u^2\Delta \zeta
    -\frac{1}{2} \io \zeta u^2\Delta v
    + \frac{1}{2} \io u^2 \nabla v \cdot \nabla \zeta \\
    &\le& \frac{1}{2} \io u^2 \Big( \Delta \zeta + \zeta w + \nabla v\cdot \nabla\zeta\Big)
    \qquad \mbox{for all } s\in (t_0,t+1).
  \eas
  Thus, by Step 4, Step 1 and (\ref{mass}), we can find $c_{12}(r_0)>0$ satisfying
  \bas
    \frac{1}{2} \frac{d}{dt} \io \zeta u^2(\cdot,s)
    + \io \zeta |\nabla u(\cdot,s)|^2
    &\le& c_{12}(r_0) \int_{\Omega\setminus B_\frac{r_0}{2}} u^2 \\
    &\le& c_{12}(r_0) \cdot m \cdot \|u(\cdot,s)\|_{L^\infty(\Omega\setminus B_\frac{r_0}{2})}
    \qquad \mbox{for all } s\in (t_0,t+1),
  \eas
  whence integrating and using (\ref{24.8}) and Step 3 shows that
  \bas
    \frac{1}{2} \io \zeta u^2(\cdot,t+1)
    + \int_{t_0}^{t+1} \io \zeta |\nabla u|^2
    &\le& \frac{1}{2} \io \zeta u^2(\cdot,t_0)
    + c_{12}(r_0) \cdot m \cdot \int_{t_0}^{t+1} \|u(\cdot,s)\|_{L^\infty(\Omega\setminus B_\frac{r_0}{2})} ds \\
    &\le& \frac{\pi}{2} c_6^2\Big(\frac{r_0}{2}\Big)
    + c_{12}(r_0) \cdot m \cdot 2c_6\Big(\frac{r_0}{2}\Big).
  \eas
  As $t_0<t$, this implies (\ref{24.7}).\abs
  \underline{Step 6.} \quad
  Conclusion.\
  Again with $\zeta$ as in Step 2, we let $\tu(r,t):=\zeta(r) u(r,t)$ for $r\in [0,1]$ and $t\ge 0$.
  Then
  \be{24.88}
    \tu_t=\tu_{rr} + f(r,t)
    \qquad \mbox{for all $r\in (0,1)$ and } t>0,
  \ee
  with
  \bas
    f(r,t):=\frac{1}{r} \zeta u_r - 2\zeta_r u_r - \zeta_{rr} u
    - \zeta u_r v_r - \mu(t) \zeta u + \zeta u w
    \qquad \mbox{for $r\in (0,1)$ and } t>0.
  \eas
  By the outcome of Step 1, Step 3, Step 4 and Step 5, for some $c_{13}(r_0)>0$ we have
  \be{24.9}
    \int_{t_0}^{t_0+2} \|f(\cdot,s)\|_{L^2((0,1))}^2 ds \le c_{13}(r_0)
    \qquad \mbox{for all } t_0\ge 1.
  \ee
  Now given $t\ge 2$, once more by Step 3 we can fix $t_0\in (t-1,t)$ fulfilling
  \be{24.10}
    \|\tu(\cdot,t_0)\|_{L^\infty((0,1))} \le c_6\Big(\frac{r_0}{2}\Big).
  \ee
  Since $\tu_r=0$ on $\partial (0,1)$, the variation-of-constants representation of $\tu$ in terms of the one-dimensional
  Neumann heat semigroup $(e^{\tau\Delta})_{\tau\ge 0}$ on the interval $(0,1)$ shows that
  \bas
    \tu(\cdot,t)=e^{(t-t_0)\Delta} \tu(\cdot,t_0)
    + \int_{t_0}^t e^{(t-s)\Delta} f(\cdot,s) ds.
  \eas
  Therefore, using standard smoothing estimates (\cite{quittner_souplet}) along with (\ref{24.10}), the H\"older inequality
  and (\ref{24.9}) we can find $c_{14}>0$ such that
  \bas
    \|\tu(\cdot,t)\|_{L^\infty((0,1))}
    &\le& \|\tu(\cdot,t_0)\|_{L^\infty((0,1))}
    + c_{14} \int_{t_0}^t (t-s)^{-\frac{1}{4}} \|f(\cdot,s)\|_{L^2((0,1))} ds \\
    &\le& c_6\Big(\frac{r_0}{2}\Big)
    + c_{14} \bigg(\int_{t_0}^t (t-s)^{-\frac{1}{2}}ds \bigg)^\frac{1}{2}
        \cdot \bigg( \int_{t_0}^t \|f(\cdot,s)\|_{L^2((0,1))}^2 ds \bigg)^\frac{1}{2} \\
    &\le& c_6\Big(\frac{r_0}{2}\Big)
    + c_{14} \cdot (2\sqrt{2})^\frac{1}{2} \cdot \big(c_{13}(r_0)\big)^\frac{1}{2}.
  \eas
  Since $\tu(r,t)=u(r,t)$ for all $r>r_0$, this establishes (\ref{24.1}).
\qed
\subsection{An $\eps$-regularity result. Proof of Theorem \ref{theo200}}\label{sect5.3}
In deriving Theorem \ref{theo200} from Lemma \ref{lem10}, we shall
make use of a regularity statement which says that solutions already
must remain bounded if only their mass concentrating in an
arbitrarily small ball centered at the origin is sufficiently small.
A first step toward this is achieved in the following lemma.
\begin{lem}\label{lem25}
  Let $\delta>0$. Then for all $p>1$ there exists $\eps=\eps(p)>0$ such that if for some $r_0\in (0,1)$,
  a radial solution of (\ref{0}) satisfies
  \be{25.1}
    \int_{B_{r_0}} u(x,t)dx < \eps
    \qquad \mbox{for all } t>0,
  \ee
  then
  \be{25.2}
    \sup_{t>0} \io u^p(x,t)dx < \infty.
  \ee
\end{lem}
\proof
  Using Young's inequality, given $p>1$ we can find $c_1=c_1(p)>0$ such that
  \be{25.3}
    \frac{p-1}{p} A^p B + AB^p \le \frac{\delta}{2} B^{p+1} + c_1 A^{p+1}
    \qquad \mbox{for all $A\ge 0$ and } B\ge 0.
  \ee
  Moreover, the Gagliardo-Nirenberg inequality says that with some $c_2=c_2(p)>0$ we have
  \be{25.4}
    \|\varphi\|_{L^\frac{2(p+1)}{p}(\Omega)}^\frac{2(p+1)}{p}
    \le c_2 \|\nabla \varphi\|_{L^2(\Omega)}^2 \|\varphi\|_{L^\frac{2}{p}(\Omega)}^\frac{2}{p}
    + c_2 \|\varphi\|_{L^\frac{2}{p}(\Omega)}^\frac{2(p+1)}{p}
    \qquad \mbox{for all } \varphi\in W^{1,2}(\Omega).
  \ee
  We claim that if (\ref{25.1}) holds with some $r_0\in (0,1)$ and
  \be{25.44}
    \eps:=\frac{p-1}{c_1 c_2 p^2},
  \ee
  then (\ref{25.2}) must be valid.
  To see this, we apply Lemma \ref{lem21} and estimate the terms on the right-hand side of (\ref{21.1}) by means of
  (\ref{25.3}) to obtain
  \bea{25.5}
    \frac{d}{dt} \bigg\{ \frac{1}{p} \io u^p + \frac{\tau}{p+1} \io w^{p+1} \bigg\}
    + \frac{4(p-1)}{p^2} \io |\nabla u^\frac{p}{2}|^2 + \delta \io w^{p+1}
    \le \frac{\delta}{2} \io w^{p+1} + c_1 \io u^{p+1}
  \eea
  for all $t>0$.
  We now fix $\zeta\in C_0^\infty(\Omega)$ such that $0\le \zeta \le 1$ in $\Omega$, $\zeta\equiv 1$ in $B_\frac{r_0}{2}$
  and $\supp \zeta \subset B_{r_0}$, and split
  \be{25.6}
    c_1\io u^{p+1} = c_1 \io \zeta^\frac{2(p+1)}{p} u^{p+1}
    + c_1 \io (1-\zeta^\frac{2(p+1)}{p}) u^{p+1},
  \ee
  where according to Lemma \ref{lem24} we can find $c_3=c_3(p,r_0)>0$ such that
  \be{25.7}
    c_1 \io (1-\zeta^\frac{2(p+1)}{p}) u^{p+1}
    \le c_1 \int_{\Omega\setminus B_\frac{r_0}{2}} u^{p+1} \le c_3
    \qquad \mbox{for all } t>0.
  \ee
  The first term on the right of (\ref{25.6}) can be estimated using (\ref{25.4}) according to
  \bea{25.8}
    c_1 \io \zeta^\frac{2(p+1)}{p} u^{p+1}
    &=& c_1 \|\zeta u^\frac{p}{2}\|_{L^\frac{2(p+1)}{p}(\Omega)}^\frac{2(p+1)}{p} \nn\\
    &\le& c_1 c_2 \|\nabla(\zeta u^\frac{p}{2})\|_{L^2(\Omega)}^2 \|\zeta u^\frac{p}{2}\|_{L^\frac{2}{p}(\Omega)}^\frac{2}{p}
    + c_1 c_2 \|\zeta u^\frac{p}{2}\|_{L^\frac{2}{p}(\Omega)}^\frac{2(p+1)}{p}.
  \eea
  Here since $\nabla (\zeta u^\frac{p}{2})=\zeta \nabla u^\frac{p}{2} + u^\frac{p}{2}\nabla \zeta$, again by Lemma \ref{lem24}
  we have
  \bas
    \|\nabla (\zeta u^\frac{p}{2})\|_{L^2(\Omega)}^2
    &\le& 2\io \zeta^2 |\nabla u^\frac{p}{2}|^2 + 2\io u^p |\nabla \zeta|^2 \\
    &\le& 2 \io |\nabla u^\frac{p}{2}|^2 + c_4
    \qquad \mbox{for all } t>0
  \eas
  with some $c_4=c_4(p,r_0)>0$, because $\supp \nabla \zeta \subset \Omega \setminus B_\frac{r_0}{2}$.
  Moreover, our hypothesis (\ref{25.1}) asserts that
  \bas
    \|\zeta u^\frac{p}{2}\|_{L^\frac{2}{p}(\Omega)}^\frac{2}{p} = \io \zeta^\frac{2}{p} u
    \le \int_{B_{r_0}} u < \eps
    \qquad \mbox{for all } t>0,
  \eas
  whence (\ref{25.8}) implies that
  \be{25.9}
    c_1 \io \zeta^\frac{2(p+1)}{p} u^{p+1} \le 2c_1 c_2 \eps \io |\nabla u^\frac{p}{2}|^2
    + c_1c_2c_4 \eps + c_1 c_2 \eps^{p+1}
    \qquad \mbox{for all } t>0.
  \ee
  Since $2c_1 c_2 \eps =\frac{2(p-1)}{p^2}$ by (\ref{25.44}), from (\ref{25.5})-(\ref{25.9}) we thus obtain
  \bea{25.10}
    \frac{d}{dt} \bigg\{ \frac{1}{p} \io u^p + \frac{\tau}{p+1} \io w^{p+1} \bigg\}
    &+& \frac{2(p-1)}{p^2} \io |\nabla u^\frac{p}{2}|^2 + \frac{\delta}{2} \io w^{p+1} \nn\\
    &\le& c_3 + c_1c_2c_4\eps + c_1 c_2 \eps^{p+1}
    \qquad \mbox{for all } t>0.
  \eea
  Here we may invoke the Poincar\'e inequality to find $c_5=c_5(p)>0$ fulfilling
  \bas
    \io \varphi^2 \le c_5 \io |\nabla \varphi|^2 + c_5 \bigg(\io |\varphi|^\frac{2}{p}\bigg)^p
    \qquad \mbox{for all } \varphi\in W^{1,2}(\Omega),
  \eas
  which according to (\ref{mass}) entails that
  \bas
    \frac{2(p-1)}{p^2} \io |\nabla u^\frac{p}{2}|^2
    \ge \frac{2(p-1)}{c_5 p^2} \io u^p - \frac{2(p-1)}{p^2} m^p
    \qquad \mbox{for all } t>0
  \eas
  with $m:=\io u_0$.  Therefore, (\ref{25.10}) shows that
  $y(t):=\frac{1}{p}\io u^p(\cdot,t)+\frac{\tau}{p+1}\io w^{p+1}$, $t\ge 0$, satisfies
  \bas
    y'(t) + c_6 y(t)\le c_7
    \qquad \mbox{for all } t>0,
  \eas
  where
  \bas
    c_6:=\min \Big\{ \frac{2(p-1)}{c_5 p} \, , \, \frac{(p+1)\delta}{2\tau} \Big\}
  \eas
  and $c_7:=c_3+c_1c_2c_4 \eps + c_1 c_2 \eps^{p+1} + \frac{2(p-1)}{p^2} m^p$.
  An ODE comparison thus yields (\ref{25.2}).
\qed
By applying the above to suitably large $p$ and using additional regularity arguments, we can next
make sure that the above assumptions already imply boundedness of $u$ with respect to the norm in
$L^\infty(\Omega)$.
\begin{lem}\label{lem26}
  Let $\delta>0$. Then there exists $\eps>0$ such that if for some $r_0\in (0,1)$ and some radial solution
  of (\ref{0}) we have
  \be{26.1}
    \int_{B_{r_0}} u(x,t)dx < \eps
    \qquad \mbox{for all } t>0,
  \ee
  then
  \be{26.2}
    \sup_{t>0} \|u(\cdot,t)\|_{L^\infty(\Omega)} < \infty.
  \ee
\end{lem}
\proof
  We pick any $p>2$ and apply Lemma \ref{lem25} which says that under the assumption (\ref{26.1}) we can find $c_1>0$ such that
  \bas
    \|u(\cdot,t)\|_{L^p(\Omega)} \le c_1
    \qquad \mbox{for all } t>0.
  \eas
  Therefore, (\ref{w_explicit}) shows that
  \bas
    \|w(\cdot,t)\|_{L^p(\Omega)}
    &\le& e^{-\frac{\delta}{\tau} t} \|w_0\|_{L^p(\Omega)}
    + \frac{1}{\tau} \int_0^t e^{-\frac{\delta}{\tau}(t-s)} \|u(\cdot,s)\|_{L^p(\Omega)} ds \\
    &\le& e^{-\frac{\delta}{\tau} t} \|w_0\|_{L^p(\Omega)}
    + \frac{c_1}{\tau} \int_0^t e^{-\frac{\delta}{\tau}(t-s)} ds
    \qquad \mbox{for all } t>0.
  \eas
  Since $\delta>0$, we know that $\int_0^t e^{-\frac{\delta}{\tau}(t-s)} ds \le \frac{\tau}{\delta}$,
  so that from this we obtain $c_2>0$ fulfilling
  \bas
    \|w(\cdot,t)\|_{L^p(\Omega)} \le c_2
    \qquad \mbox{for all } t>0.
  \eas
  From this and standard
  elliptic regularity theory we obtain a bound for $v$ in
  $L^\infty((0, \infty); W^{2, p}(\Omega))$, which by the validity of the embedding
  $W^{2,p}(\Omega) \hookrightarrow W^{1,\infty}(\Omega)$ implies that
  \bas
    \|\nabla v(\cdot, t)\|_{L^\infty(\Omega)}\le c_3 \qquad \mbox{for all $t\in (0, T)$}
  \eas
  with some $c_3>0$. Combined with (\ref{25.2}), upon a Moser-type iteration (\cite[Lemma 4.1]{taowin2}) this
  yields (\ref{26.2}).
\qed
Combining Lemma \ref{lem26} with Lemma \ref{lem10} now immediately yields boundedness of solutions in the subcritical
mass case.\abs
\proofc of Theorem \ref{theo200}. \quad
  We let $\eps>0$ be as provided by Lemma \ref{lem26} and only need to verify the validity of (\ref{26.1})
  for some $r_0\in (0,1)$.
  In order to choose the latter appropriately, we apply Lemma \ref{lem10} to find $c_1>0$ such that for arbitrary
  $r_0\in (0,1)$ we have
  \bas
    \int_{B_{r_0}} u(x,t)dx = 2\pi \int_0^{r_0} r(u(r, t) dr =2\pi
    U(r_0^2, t)\le c_1 r_0^2 \qquad\mbox{for all $t>0$}.
  \eas
  This means that if we now fix $r_0\in (0,1)$ in such a way that $r_0<\sqrt{\frac{\eps}{c_1}}$, then indeed
  \bas
    \int_{B_{r_0}} u(x,t)dx <\eps \qquad\mbox{for all $t>0$}.
  \eas
  Lemma \ref{lem26} thus ensures that (\ref{26.2}) holds, whereupon recalling (\ref{w_explicit}) and applying
  elliptic regularity theory we see that the statement in Theorem \ref{theo200} becomes an evident consequence thereof.
\qed
\mysection{Unbounded solutions with $\io u_0>8\pi\delta$. Proof of Theorem \ref{theo300}}\label{sect6}
\subsection{A class of comparison functions}
We shall next prove that whenever $m>8\pi\delta$, some solutions at the mass level $m$
asymptotically aggregate in the spirit of Theorem \ref{theo300}.
To this end, we shall consider comparison functions $\UU:[0,1]\times [0,\infty)\to\R$ of the form
\be{UU}
    \UU(\xi,t):=\left\{ \begin{array}{ll}
    \frac{a(t)\xi}{b(t)+\xi} & \mbox{if } \xi\in [0,\xi_0] \mbox{ and } t\ge 0, \\[2mm]
    \frac{a(t)b(t)\xi+a(t)\xi_0^2}{(b(t)+\xi_0)^2} \qquad & \mbox{if $\xi\in (\xi_0,1]$ and } t\ge 0,
    \end{array} \right.
\ee
where $\xi_0\in (0,1)$ and $a$ and $b$ are a suitably chosen positive functions on $[0,\infty)$.
Let us first collect some basic properties of such functions, especially with regard to their behavior under
the action of the operator $\parab$ defined in (\ref{parab}).
\begin{lem}\label{lem5}
  Let $\xi_0\in (0,1)$, and assume that $a\in C^1([0,\infty))$ and $b\in C^1([0,\infty))$ are positive.
  Then the function $\UU$ given by (\ref{UU}) satisfies
  \bas
    \UU \in C^1([0,1]\times [0,\infty)) \cap C^0([0,\infty);W^{2,\infty}((0,1)))
    \cap C^0([0,\infty); C^2_{loc}([0,1]\setminus \{\xi_0\})).
  \eas
  Moreover, with $\parab$ as in (\ref{parab}) we have
  \bea{5.2}
    \frac{(b(t)+\xi)^2}{a(t)b(t)\xi} \cdot \parab \UU(\xi,t)
    &=& \frac{a'(t) (b(t)+\xi)}{a(t)b(t)}
    -\frac{b'(t)}{b(t)} + \frac{8}{b(t)+\xi} \nn\\
    & & - \frac{2}{\tau}\int_0^t e^{-\frac{\delta}{\tau}(t-s)} \Big\{ \frac{a(s)}{b(s)+\xi} - \frac{m}{2\pi}\Big\} ds \nn\\
    & & - 2\Big(\frac{W_0(\xi)}{\xi}-\kz\Big)\cdot e^{-\frac{\delta}{\tau} t}
    \qquad \mbox{for all $\xi\in (0,\xi_0)$ and $t>0$}
  \eea
  and
  \bea{5.3}
    \frac{(b(t)+\xi_0)^2}{a(t)b(t)} \cdot \parab \UU(\xi,t)
    &=& \frac{a'(t)\xi}{a(t)} + \frac{b'(t)\xi}{b(t)} + \frac{a'(t) \xi_0^2}{a(t)b(t)}
    - 2\frac{b'(t)\xi + \frac{b'(t)}{b(t)} \xi_0^2}{b(t)+\xi_0} \nn\\
    & & - \frac{2}{\tau} \int_0^t e^{-\frac{\delta}{\tau}(t-s)} \bigg\{ \frac{a(s) b(s)\xi+a(s) \xi_0^2}{(b(s)+\xi_0)^2}
        - \frac{m}{2\pi}\xi\bigg\} ds \nn\\
    & & - 2\Big(W_0(\xi)-\kz \xi\Big)\cdot e^{-\frac{\delta}{\tau} t}
    \qquad \mbox{for all $\xi\in (\xi_0,1)$ and $t>0$,}
  \eea
  where $W_0$ and $\kz$ are as defined in (\ref{W0}).
\end{lem}
\proof
  The claimed regularity properties can immediately be verified using the explicit form of $\UU$ which
  clearly allows for piecewise differentiation, resulting in
  \be{5.4}
    U_\xi(\xi,t)=\left\{ \begin{array}{ll}
    \frac{ab}{(b+\xi)^2} \qquad & \mbox{for $\xi\in (0,\xi_0)$ and } t>0, \\[2mm]
    \frac{ab}{(b+\xi_0)^2} \qquad & \mbox{for $\xi\in (\xi_0,1)$ and } t>0,
    \end{array} \right.
  \ee
  and
  \be{5.5}
    U_{\xi\xi}(\xi,t)=\left\{ \begin{array}{ll}
    - \frac{2ab}{(b+\xi)^3} \qquad & \mbox{for $\xi\in (0,\xi_0)$ and } t>0, \\[2mm]
    0 \qquad & \mbox{for $\xi\in (\xi_0,1)$ and } t>0,
    \end{array} \right.
  \ee
  as well as
  \be{5.6}
    U_t(\xi,t)=\left\{ \begin{array}{ll}
    \frac{a'\xi}{b+\xi} - \frac{ab'\xi}{(b+\xi)^2} \qquad & \mbox{for $\xi\in (0,\xi_0)$ and } t>0, \\[2mm]
    \frac{a'b\xi + ab'\xi + a'\xi_0^2}{(b+\xi_0)^2} - 2\frac{abb'\xi+ab'\xi_0^2}{(b+\xi_0)^3}
    \qquad & \mbox{for $\xi\in (\xi_0,1)$ and } t>0.
    \end{array} \right.
  \ee
  Moreover, for $\xi<\xi_0$ we obtain from (\ref{5.4})-(\ref{5.6}) that
  \bas
    \parab \UU(\xi,t)
    &=& \frac{a'\xi}{b+\xi} - \frac{ab'\xi}{(b+\xi)^2} + \frac{8ab\xi}{(b+\xi)^3} \\
    & & -\frac{2}{\tau}\Bigg\{ \int_0^t e^{-\frac{\delta}{\tau}(t-s)} \cdot
    \Big\{ \frac{a(s)\xi}{b(s)+\xi} - \frac{m}{2\pi}\xi \Big\} ds \Bigg\}
        \cdot \frac{ab}{(b+\xi)^2} \\
    & & - \Bigg\{ 2\Big(W_0(\xi)-\kz \xi\Big)\cdot e^{-\frac{\delta}{\tau} t} \Bigg\} \cdot \frac{ab}{(b+\xi)^2},
  \eas
  which is equivalent to (\ref{5.2}).
  Likewise, (\ref{5.3}) easily follows upon the observation that for $\xi>\xi_0$, the identity
  \bas
    \frac{(b+\xi_0)^2}{ab} \cdot \Big(\UU_t-4\xi\UU_{\xi\xi}\Big)
    &=& \frac{a'b\xi+ab'\xi+a'\xi_0^2}{ab} -2\frac{abb'\xi +ab'\xi_0^2}{ab(b+\xi_0)} \\
    &=& \frac{a'\xi}{a} + \frac{b'\xi}{b} + \frac{a'\xi_0^2}{ab}
    -2\frac{b'\xi+\frac{b'}{b} \xi_0^2}{b+\xi_0}
  \eas
  holds.
\qed
To make the above choice as efficient as possible for our purpose, $a(t)$ will be adjusted in such a way that
at $\xi=1$, the function $\UU$ attains the same boundary value as $U$ introduced in (\ref{31.1}).
The corresponding condition $\UU(1,t)=\frac{m}{2\pi}$ for all $t\ge 0$, with $m:=\io u_0$,
thus amounts to requiring that $a(t)$
is linked to $\xi_0$ and $b(t)$, and accordingly we shall concentrate on the case when
\be{a}
    a(t):=\frac{m}{2\pi} \cdot \frac{(b(t)+\xi_0)^2}{b(t)+\xi_0^2}
    \qquad \mbox{for } t\ge 0
\ee
in the sequel.
Then for later use we note that if in addition we assume that $b$ is differentiable, so will be $a$ with
\bea{aprime}
    a'(t)
    &=& \frac{m}{2\pi} \cdot \frac{2(b+\xi_0)(b+\xi_0^2) - (b+\xi_0)^2}{(b+\xi_0^2)^2} \cdot b' \nn\\
    &=& \frac{m}{2\pi} \cdot \frac{b^2+2b\xi_0^2-\xi_0^2 + 2\xi_0^3}{(b+\xi_0^2)^2} \cdot b'
    \qquad \mbox{for all } t>0,
\eea
\subsection{Subsolution in an annulus}
We first analyze in more depth the behavior of $\UU$ in the outer region where $\xi>\xi_0$.
Here it will not be necessary to fix $\xi_0$, and keeping this freedom will be important for our procedure
in the corresponding inner part where $\xi<\xi_0$, in which we shall adjust $\xi_0$ in dependence of $m>8\pi\delta$.\\
To begin with, let us draw a first conclusion from Lemma \ref{lem5} under the assumption (\ref{a}).
\begin{lem}\label{lem6}
  Let $\delta\ge 0$ and $m>0$, and
  suppose that $\xi_0\in (0,1)$, that $b\in C^1([0,\infty))$ is positive and nonincreasing such that
  \be{6.1}
    b(t)\le \xi_0^2
    \qquad \mbox{for all } t\ge 0,
  \ee
  and that $a\in C^1([0,\infty))$ is given by (\ref{a}).
  Then the function $\UU$ defined in (\ref{UU}) satisfies
  \bea{6.3}
    \frac{(b(t)+\xi_0)^2}{a(t)b(t)} \cdot \parab \UU(\xi,t)
    &\le& (1-\xi)\cdot \bigg\{ -\frac{b'(t)}{b(t)} - \frac{m}{2\pi\tau} \int_0^t e^{-\frac{\delta}{\tau}(t-s)} ds \bigg\}
    - 2\Big(W_0(\xi)-\kz \xi \Big) \cdot e^{-\frac{\delta}{\tau} t} \nn\\[2mm]
    & & \hspace*{40mm}  \mbox{for all $\xi\in (\xi_0,1)$ and } t>0,
  \eea
  where $\parab$ is as in (\ref{parab}).
\end{lem}
\proof
  From (\ref{aprime}) we compute
  \bas
    \frac{a'(t)}{a(t)}
    = \frac{b^2+2b\xi_0^2-\xi_0^2 + 2\xi_0^3}{(b+\xi_0^2)(b+\xi_0)^2} \cdot b'
    \qquad \mbox{for all } t>0.
  \eas
  Thus, on the right-hand side of (\ref{5.3}) we have
  \bas
    J_1(\xi,t)
    &:=& \frac{a' \xi}{a} + \frac{b'\xi}{b} + \frac{a'\xi_0^2}{ab} - 2\frac{b'\xi + \frac{b'}{b}\xi_0^2}{b+\xi_0} \\
    &=& \frac{b'}{b} \cdot \Bigg\{ \frac{(b^2+2b\xi_0^2 - \xi_0^2 + 2\xi_0^3)(b\xi+\xi_0^2)}{(b+\xi_0^2)(b+\xi_0)^2}
    + \xi - \frac{2b\xi+2\xi_0^2}{b+\xi_0} \Bigg\}
    \qquad \mbox{for all $\xi\in (\xi_0,1)$ and } t>0,
  \eas
  whereupon a lengthy but straightforward computation yields
  \be{6.4}
    J_1(\xi,t)
    = - \frac{b'}{b} \cdot \frac{\xi_0^2 (1-\xi)}{b+\xi_0^2}
    \qquad \mbox{for all $\xi\in (\xi_0,1)$ and } t>0.
  \ee
  Next, in the integrand on the right of (\ref{5.3}) we again use (\ref{a}) to see that
  \bas
    \frac{ab\xi+a\xi_0^2}{(b+\xi_0)^2} - \frac{m}{2\pi} \xi
    &=& \frac{m}{2\pi} \cdot \bigg( \frac{b\xi+\xi_0^2}{b+\xi_0^2}-\xi\bigg) \\
    &=& \frac{m}{2\pi} \cdot \frac{\xi_0^2(1-\xi)}{b+\xi_0^2}
    \qquad \mbox{for all $\xi\in (\xi_0,1)$ and } t>0,
  \eas
  so that
  \bea{6.5}
    J_2(\xi,t)
    &:=& -\frac{2}{\tau} \int_0^t e^{-\frac{\delta}{\tau}(t-s)} \cdot
    \bigg( \frac{a(s)b(s)\xi+a(s)\xi_0^2}{(b(s)+\xi_0)^2} - \frac{m}{2\pi}\xi\bigg) ds \nn\\
    &=& -\frac{m}{\pi\tau} \cdot (1-\xi)\cdot \int_0^t e^{-\frac{\delta}{\tau}(t-s)} \cdot \frac{\xi_0^2}{b(s)+\xi_0^2} ds
    \qquad \mbox{for all $\xi\in (\xi_0,1)$ and } t>0.
  \eea
  Now in (\ref{6.4}) we can use the nonnegativity of $b$ to find that
  \bas
    \frac{\xi_0^2(1-\xi)}{b+\xi_0^2} \le 1-\xi
    \qquad \mbox{for all $\xi\in (\xi_0,1)$ and } t>0,
  \eas
  whereas in (\ref{6.5}) we employ (\ref{6.1}) to estimate
  \bas
    \frac{\xi_0^2}{b+\xi_0^2} \ge \frac{1}{2}
    \qquad \mbox{for all $\xi\in (\xi_0,1)$ and } t>0.
  \eas
  Therefore, by means of the nonpositivity of $\frac{b'}{b}$ we have
  \bas
    J_1(\xi,t)+J_2(\xi,t)
    \le -\frac{b'}{b} \cdot (1-\xi)
    -\frac{m}{2\pi\tau} \cdot (1-\xi) \cdot \int_0^t e^{-\frac{\delta}{\tau}(t-s)} ds
    \qquad \mbox{for all $\xi\in (\xi_0,1)$ and } t>0.
  \eas
  In view of (\ref{5.3}), this proves (\ref{6.3}).
\qed
Now the right-hand side of (\ref{6.3}) suggests to choose $b$ in such a way that $\frac{b'}{b}$ is a negative constant.
In that case, namely, it turns out that
the unfavorable contribution of $-\frac{b'}{b}$ in (\ref{6.3}) can be controlled for large times by the
integral on the right of (\ref{6.3}), whereas for small $t$ it will be dominated by the expression containing
$W_0$ and $\kz$, provided that $w_0$ satisfies some rather mild condition.
\begin{lem}\label{lem2}
  Let $\delta\ge 0$ and $m> 0$, and suppose that for some $\xi_0\in (0,1)$ and $\eta_0>0$, with $W_0$ and $\kz$ as in
  (\ref{W0}) we have
  \be{2.1}
    \frac{W_0(\xi)-\kz\xi}{1-\xi} \ge \eta_0
    \qquad \mbox{for all } \xi \in (\xi_0,1).
  \ee
  Then for all $\alpha_\star>0$ there exists $\alpha \in (0,\alpha_\star)$
  such that for any choice of $b_0\in (0,\xi_0^2)$, with
  \bas
    b(t):=b_0 e^{-\alpha t}
    \qquad \mbox{for } t\ge 0
  \eas
  and $a\in C^1([0,\infty))$ as in (\ref{a}), the function $\UU$ in (\ref{UU}) satisfies
  \be{2.2}
    \parab \UU(\xi,t) \le 0
    \qquad \mbox{for all $\xi\in (\xi_0,1)$ and } t>0,
  \ee
  the operator $\parab$ being defined through (\ref{parab}).
\end{lem}
\proof
  We claim that (\ref{2.2}) holds whenever $b_0\in (0,\xi_0^2)$ and
  \be{2.4}
    \alpha<\alpha_0 :=\min \bigg\{ \frac{m}{2\pi\tau e^\frac{\delta}{\tau}} \, , \,
    \frac{2\eta_0}{e^\frac{2\delta}{\tau}} \, , \, \alpha_\star \bigg\} \, .
  \ee
  Indeed, since (\ref{2.1}) in particular implies that $W_0(\xi)-\kz \xi \ge 0$ for all $\xi\in (\xi_0,1)$,
  from Lemma \ref{lem6} we obtain that
  \bas
    \frac{(b(t)+\xi_0)^2}{a(t)b(t)} \cdot \parab \UU(\xi,t)
    &\le& (1-\xi) \cdot \bigg\{ - \frac{b'(t)}{b(t)} -\frac{m}{2\pi\tau} \int_0^t e^{-\frac{\delta}{\tau}(t-s)} ds \bigg\} \\
    &=& (1-\xi) \cdot \bigg\{ \alpha - \frac{m}{2\pi\tau} \int_0^t e^{-\frac{\delta}{\tau}(t-s)} ds \bigg\}
    \qquad \mbox{for all $\xi\in (\xi_0,1)$ and } t>0.
  \eas
  Here for large $t$ we can estimate
  \bas
    \int_0^t e^{-\frac{\delta}{\tau}(t-s)} ds
    \ge \int_{t-1}^t e^{-\frac{\delta}{\tau}(t-s)} ds
    \ge \int_0^1 e^{-\frac{\delta}{\tau}} ds = e^{-\frac{\delta}{\tau}}
    \qquad \mbox{for all } t\ge 2,
  \eas
  so that the first restriction implied by (\ref{2.4}) warrants that
  \bas
    \frac{(b(t)+\xi_0)^2}{a(t)b(t)} \cdot \parab \UU(\xi,t)
    &\le& (1-\xi) \cdot \bigg\{ \alpha - \frac{m}{2\pi\tau} e^{-\frac{\delta}{\tau}} \bigg\} \\[2mm]
    &\le& 0
    \qquad \mbox{for all $\xi\in (\xi_0,1)$ and } t\ge 2.
  \eas
  For small values of $t$, however, (\ref{6.3}) and (\ref{2.1}) yield
  \bas
    \frac{(b(t)+\xi_0)^2}{a(t)b(t)} \cdot \parab \UU(\xi,t)
    &\le& (1-\xi)\alpha - 2\eta_0 (1-\xi) e^{-\frac{\delta}{\tau} t} \\
    &\le& (1-\xi)\alpha - 2\eta_0(1-\xi) e^{-2\frac{\delta}{\tau}} \\[2mm]
    &\le& 0
    \qquad \mbox{for all $\xi\in (\xi_0,1)$ and } t< 2
  \eas
  because of the second limitation on $\alpha$ asserted by (\ref{2.4}).
\qed
\subsection{Subsolution near the origin}
Our argument in the associated inner region will be more subtle, and here we will in particular rely on the supercriticality
assumption $m>8\pi\delta$.
Let us begin by estimating the first term on the right of (\ref{5.2}) under the hypothesis (\ref{a}).
\begin{lem}\label{lem3}
  Let $m> 0$, and suppose that $b \in C^1([0,\infty)$ is positive and nonincreasing, and let $\xi_0\in (0,1)$.
  Then the function $a\in C^1([0,\infty))$ defined in (\ref{a}) satisfies
  \bas
    \frac{a'(t)(b(t)+\xi)}{a(t)b(t)} \le \frac{1}{\xi_0} \cdot \frac{-b'(t)}{b(t)}
    \qquad \mbox{for all $\xi\in (0,\xi_0)$ and } t>0.
  \eas
\end{lem}
\proof
  In (\ref{aprime}) we can trivially estimate
  \bas
    b^2(t)+2b(t)\xi_0^2 - \xi_0^2 + \xi_0^3 \ge -\xi_0^2
    \qquad \mbox{for all } t>0
  \eas
  to obtain
  \bas
    a'(t) \le - \frac{m}{2\pi} \cdot \frac{\xi_0^2}{(b(t)+\xi_0^2)^2} \cdot b'(t)
    \qquad \mbox{for all } t>0.
  \eas
  Therefore
  \bas
    \frac{a'(t)(b(t)+\xi)}{a(t)b(t)}
    &\le& \frac{-\frac{\xi_0^2}{(b(t)+\xi_0^2)^2} \cdot b'(t)}{\frac{(b(t)+\xi_0)^2}{b(t)+\xi_0^2}}
        \cdot \frac{b(t)+\xi}{b(t)} \\
    &=& \frac{\xi_0^2 (b(t)+\xi)}{(b(t)+\xi_0^2)(b(t)+\xi_0)^2} \cdot \frac{-b'(t)}{b(t)}
    \qquad \mbox{for all $\xi\in (0,1)$ and } t>0,
  \eas
  so that since $\frac{b(t)+\xi}{b(t)+\xi_0} \le 1$ for all $\xi\in (0,\xi_0)$ and $t>0$, we find that
  \bas
    \frac{a'(t)(b(t)+\xi)}{a(t)b(t)}
    &\le& \frac{\xi_0^2}{(b(t)+\xi_0^2)(b(t)+\xi_0)} \cdot \frac{-b'(t)}{b(t)} \\
    &\le& \frac{1}{\xi_0} \cdot \frac{-b'(t)}{b(t)}
    \qquad \mbox{for all $\xi\in (0,\xi_0)$ and } t>0,
  \eas
  again because $b\ge 0$ and $b'\le 0$.
\qed
The technical key toward our proof of infinite-time blow-up in
the supercritical case is contained in the following lemma which
says that in the supercritical mass case we can achieve that $\UU$
is a subsolution in the inner region for suitably large times upon
an appropriate choice of the parameters.
\begin{lem}\label{lem1}
  Let $\delta\ge 0$ and
  \be{1.1}
    m>8\pi\delta,
  \ee
  and suppose that taking $W_0$ and $\kz$ from (\ref{W0}), we have
  \be{1.01}
    W_0(\xi)-\kz \xi \ge 0
    \qquad \mbox{for all } \xi\in (0,1).
  \ee
  Then there exist $\xi_0\in (0,1)$ and $\alpha_\star>0$ with the property that for all $\alpha\in (0,\alpha_\star)$
  one can find $b_0\in (0,\xi_0^2)$ and $t_0>0$ such that with
  \be{b_exp}
    b(t):=b_0 e^{-\alpha t}
    \qquad \mbox{for } t\ge 0
  \ee
  and $a\in C^1([0,\infty))$ as given by (\ref{a}), the function $\UU$ in (\ref{UU}) satisfies
  \be{1.2}
    \parab \UU(\xi,t) \le 0
    \qquad \mbox{for all $\xi\in (0,\xi_0)$ and } t\ge t_0,
  \ee
  where $\parab$ is given by (\ref{parab}).
\end{lem}
\proof
  We detail the proof for the case when $\delta$ is positive, leaving the minor modifications necessary for
  the limit case $\delta=0$ to the reader.
  Since $m>8\pi\delta$, we can pick $\eps\in (0,1)$ small enough such that
  \be{1.3}
    c_1:=\frac{(1-\eps)^3 m}{(1+\eps)\pi\delta} - 8 >0,
  \ee
  and thereafter fix $\alpha_\star>0$ and $\xi_0\in (0,1)$ small fulfilling
  \be{1.4}
    \alpha_\star \le \frac{c_1}{4}
  \ee
  and
  \be{1.5}
    \alpha_\star \le \frac{\delta \ln \frac{1}{1-\eps}}{\tau\ln \frac{1}{\eps}}
  \ee
  as well as
  \be{1.6}
    \xi_0 \le \frac{\eps}{2}.
  \ee
  Given $\alpha \in (0,\alpha_\star)$, we then choose $b_0>0$ suitably small and $t_0>0$ sufficiently large such that
  \be{1.7}
    b_0 \le \eps \xi_0^2
  \ee
  and
  \be{1.8}
    t_0 \ge \frac{1}{\alpha} \cdot \ln \frac{1}{1-\eps},
  \ee
  and thereupon let $b, a$ and $\UU$ be defined by (\ref{b_exp}), (\ref{a}) and (\ref{UU}).
  Then by (\ref{1.01}), Lemma \ref{lem5} implies that
  \bea{1.88}
    \frac{(b(t)+\xi)^2}{a(t)b(t)\xi} \cdot \parab \UU(\xi,t)
    \le J(\xi,t)
    := \frac{a'(b+\xi)}{ab} - \frac{b'}{b} + \frac{8}{b+\xi}
    - \frac{2}{\tau}\int_0^t e^{-\frac{\delta}{\tau}(t-s)} \cdot \bigg\{ \frac{a(s)}{b(s)+\xi} - \frac{m}{2\pi} \bigg\} ds
  \eea
  for all $\xi\in (0,\xi_0)$ and $t>0$.
  Here by Lemma \ref{lem3}, we can estimate
  \bea{1.9}
    J_1(\xi,t)
    &:=& \frac{a'(b+\xi)}{ab} - \frac{b'}{b} \nn\\
    &\le& - \bigg(\frac{1}{\xi_0}+1\bigg) \cdot \frac{b'}{b}  \nn\\
    &=&  \bigg(\frac{1}{\xi_0}+1\bigg) \cdot \alpha \nn\\
    &\le& \frac{2}{\xi_0}\cdot \alpha
    \qquad \mbox{for all } t>0.
  \eea
  Next, to estimate the integral in (\ref{1.88}) we first note that (\ref{1.6}) guarantees that
  \bas
    \eps \cdot \frac{a(t)}{b(t)+\xi} - \frac{m}{2\pi}
    &=& \frac{m}{2\pi} \cdot \bigg\{ \eps \cdot \frac{b+\xi_0}{b+\xi} \cdot \frac{b+\xi_0}{b+\xi_0^2} - 1 \bigg\} \\
    &\ge& \frac{m}{2\pi} \cdot \bigg\{ \eps\cdot 1 \cdot \frac{\xi_0}{2\xi_0^2} - 1 \bigg\} \\[2mm]
    &\ge& 0
    \qquad \mbox{for all $\xi  \in (0,\xi_0)$ and } t>0,
  \eas
  and that
  \bas
    a(t) \ge \frac{m}{2\pi} \cdot \frac{\xi_0^2}{b+\xi_0^2}
    \ge \frac{m}{2\pi} \cdot \frac{\xi_0^2}{(1+\eps)\xi_0^2}
    = \frac{m}{2(1+\eps)\pi}
    \qquad \mbox{for all } t>0
  \eas
  by (\ref{1.7}). Therefore,
  \bea{1.10}
    J_2(\xi,t)
    &:=& \frac{8}{b(t)+\xi} -
    \frac{2}{\tau}\int_0^t e^{-\frac{\delta}{\tau}(t-s)} \cdot \bigg\{ \frac{a(s)}{b(s)+\xi} - \frac{m}{2\pi} \bigg\} ds
        \nn\\
    &\le& \frac{8}{b(t)+\xi}
    - \frac{2(1-\eps)}{\tau} \int_0^t e^{-\frac{\delta}{\tau}(t-s)} \cdot \frac{a(s)}{b(s)+\xi} ds \nn\\
    &\le& \frac{8}{b(t)+\xi}
    - \frac{(1-\eps)m}{(1+\eps)\pi\tau} \int_0^t e^{-\frac{\delta}{\tau}(t-s)} \cdot \frac{1}{b(s)+\xi} ds
    \qquad \mbox{for all $\xi  \in (0,\xi_0)$ and } t>0.
  \eea
  Now whenever $t\ge s \ge t-\frac{1}{\alpha} \ln \frac{1}{1-\eps} \ge 0$, we have
  \bas
    \frac{b(t)}{b(s)} = e^{-\alpha(t-s)} \ge e^{-\ln\frac{1}{1-\eps}} =1-\eps,
  \eas
  which implies that
  \bas
    \frac{b(t)+\xi}{b(s)+\xi}
    \ge \frac{b(t)+(1-\eps)\xi}{b(s)+\xi} \ge 1-\eps
  \eas
  for all $\xi>0$ and any such $t$ and $s$.
  By means of (\ref{1.8}), we can hence estimate
  \bas
    \int_0^t e^{-\frac{\delta}{\tau}(t-s)} \cdot \frac{1}{b(s)+\xi} ds
    &\ge& \frac{1-\eps}{b(t)+\xi} \cdot \int_{t-\frac{1}{\alpha} \ln \frac{1}{1-\eps}}^t e^{-\frac{\delta}{\tau}(t-s)} ds \\
    &=& \frac{1-\eps}{b(t)+\xi} \cdot \frac{\tau}{\delta}
    \bigg( 1 - e^{-\frac{\delta}{\alpha\tau} \ln \frac{1}{1-\eps}} \bigg) \\
    &\ge& \frac{(1-\eps)^2 \tau}{(b(t)+\xi) \cdot\delta}
    \qquad \mbox{for all $\xi  \in (0,\xi_0)$ and } t\ge t_0,
  \eas
  because (\ref{1.5}) ensures that
  \bas
    e^{-\frac{\delta}{\alpha\tau} \ln \frac{1}{1-\eps}}
    \le \exp \bigg( -\frac{\delta}{\tau}\ln\frac{1}{1-\eps}
    \cdot \frac{\tau \ln\frac{1}{\eps}}{\delta \ln\frac{1}{1-\eps}} \bigg) = \eps.
  \eas
  Accordingly, (\ref{1.10}) and (\ref{1.3}) yield
  \bas
    J_2(\xi,t)
    &\le& \frac{1}{b(t)+\xi} \cdot
    \bigg\{ 8 - \frac{(1-\eps)m}{(1+\eps)\pi\tau} \cdot \frac{(1-\eps)^2 \tau}{\delta} \bigg\} \\
    &=& - \frac{c_1}{b(t)+\xi}
    \qquad \mbox{for all $\xi  \in (0,\xi_0)$ and } t\ge t_0,
  \eas
  which since by (\ref{1.7}) we have
  \bas
    b(t)+\xi \le b_0+\xi_0 \le \eps \xi_0^2 + \xi_0 \le 2\xi_0
    \qquad \mbox{for all $\xi  \in (0,\xi_0)$ and } t>0
  \eas
  guarantees that
  \bas
    J_2(\xi,t)
    \le - \frac{c_1}{2\xi_0}
    \qquad \mbox{for all $\xi \in (0,\xi_0)$ and } t\ge t_0.
   \eas
  In conjunction with (\ref{1.9}), (\ref{1.88}) and (\ref{1.4}), this shows that
  \bas
    J(\xi,t)
    \le \frac{2}{\xi_0} \cdot \alpha - \frac{c_1}{2\xi_0} \le 0
    \qquad \mbox{for all $\xi \in (0,\xi_0)$ and } t\ge t_0
  \eas
  and thereby proves (\ref{1.2}).
\qed
Now for small times, we can also achieve that $\parab \UU \le 0$ if we suppose $w_0$ to be sufficiently strongly
concentrated near the origin.
\begin{lem}\label{lem4}
  Let $\delta\ge 0$ and $m>0$, and suppose that $\alpha>0, b_0>0$ and $\xi_0\in (0,1)$, that
  \be{b}
    b(t)=b_0 e^{-\alpha t}
    \qquad \mbox{for } t\ge 0,
  \ee
  and that $a\in C^1([0,\infty))$ is as given by (\ref{a}).
  Then for all $t_0>0$ there exists $\Gamma_0(\alpha,b_0,\xi_0,t_0)>0$
  such that whenever $W_0$ and $\kz$ as introduced in (\ref{W0}) satisfy
  \be{4.1}
    \frac{W_0(\xi)}{\xi}-\kz \ge \Gamma_0(\alpha,b_0,\xi_0,t_0)
    \qquad \mbox{for all } \xi\in (0,\xi_0),
  \ee
  the function $\UU$ defined in (\ref{UU}) has the property that with $\parab$ as in (\ref{parab}) we have
  \be{4.2}
    \parab \UU(\xi,t) \le 0
    \qquad \mbox{for all $\xi\in (0,\xi_0)$ and $t\in (0,t_0)$.}
  \ee
\end{lem}
\proof
  By (\ref{5.2}), we have
  \bea{4.3}
    \frac{(b+\xi)^2}{ab\xi} \cdot \parab \UU(\xi,t)
    &=& \frac{a'(b+\xi)}{ab} - \frac{b'}{b} + \frac{8}{b+\xi}
    - \frac{2}{\tau}\int_0^t e^{-\frac{\delta}{\tau}(t-s)} \cdot
    \bigg\{ \frac{a(s)}{b(s)+\xi} - \frac{m}{2\pi}\bigg\} ds \nn\\
    & & - 2\bigg(\frac{W_0(\xi)}{\xi}-\kz\bigg)\cdot e^{-\frac{\delta}{\tau} t}
    \qquad \mbox{for all $\xi\in (0,\xi_0)$ and } t>0,
  \eea
  where Lemma \ref{lem3} ensures that
  \be{4.4}
    \frac{a'(b+\xi)}{ab} - \frac{b'}{b}
    \le \bigg(\frac{1}{\xi_0}+1\bigg)\cdot \frac{-b'}{b}
    = \bigg(\frac{1}{\xi_0}+1\bigg)\cdot \alpha,
  \ee
  and where
  \be{4.44}
    \frac{8}{b+\xi} \le \frac{8}{b_0} e^{\alpha t_0}
    \qquad \mbox{for all $\xi\in (0,\xi_0)$ and } t\in (0,t_0)
  \ee
  according to (\ref{b}).
  Moreover, due to (\ref{a}) we see that
  \be{4.5}
    \frac{a}{b+\xi} - \frac{m}{2\pi}
    = \frac{m}{2\pi} \cdot \frac{b+\xi_0}{b+\xi} \cdot \frac{b+\xi_0}{b+\xi_0^2} - \frac{m}{2\pi}
    \ge 0
    \qquad \mbox{for all $\xi\in (0,\xi_0)$ and } t>0,
  \ee
  because for any choice of $\xi<\xi_0$ we have $b+\xi_0\ge b+\xi$ and also $b+\xi_0\ge b+\xi_0^2$ due to the fact
  that $\xi_0<1$. Therefore,
  \bas
    - \frac{2}{\tau}\int_0^t e^{-\frac{\delta}{\tau}(t-s)}
    \cdot \bigg\{ \frac{a(s)}{b(s)+\xi} - \frac{m}{2\pi}\bigg\} ds \le 0
    \qquad \mbox{for all $\xi\in (0,\xi_0)$ and } t>0,
  \eas
  which combined with (\ref{4.3})-(\ref{4.5}) shows that
  \be{4.6}
    \frac{ab\xi}{(b+\xi)^2} \cdot \parab \UU(\xi,t)
    \le \bigg(\frac{1}{\xi_0}+1\bigg)\cdot\alpha + \frac{8}{b_0} e^{\alpha t_0}
    - 2 \bigg( \frac{W_0(\xi)}{\xi}-\kz\bigg) \cdot e^{-\frac{\delta}{\tau} t}
    \qquad \mbox{for all $\xi\in (0,\xi_0)$ and } t\in (0,t_0).
  \ee
  Thus, if
  \bas
    \frac{W_0(\xi)}{\xi}-\kz
    \ge \Gamma_0(\alpha,b_0,\xi_0,t_0)
    :=\frac{1}{2} \cdot \bigg\{ \Big(\frac{1}{\xi_0}+1\Big)\cdot\alpha + \frac{8}{b_0} e^{\alpha t_0} \bigg\}
        \cdot e^{\frac{\delta}{\tau} t_0}
    \qquad \mbox{for all } \xi \in (0,\xi_0),
  \eas
  then (\ref{4.2}) results from (\ref{4.6}).
\qed
By a careful selection of the parameters in (\ref{UU}) we can finally combine Lemma \ref{lem2}, Lemma \ref{lem1}
and Lemma \ref{lem4} to establish our main result on infinite-time blow-up of supercritical-mass solutions.\abs
\proofc of Theorem \ref{theo300}. \quad
  We first take $\xi_0 \in (0,1)$ and $\alpha_\star>0$ as provided by Lemma \ref{lem1} and let
  \be{300.4}
    R:=\sqrt{\xi_0}.
  \ee
  Again invoking Lemma \ref{lem1}, we then find $b_0\in (0,\xi_0^2)$ and $t_0>0$ with the properties listed there, and
  thereupon pick $\alpha \in (0,\alpha_\star)$ as given by Lemma \ref{lem2} when applied to $\eta_0:=\frac{\eta}{2}$.
  Having thus fixed $\alpha, b_0, \xi_0$ and $t_0$, we finally fix $\Gamma_0:=\Gamma_0(\alpha,b_0,\xi_0,t_0)$ as yielded by
  Lemma \ref{lem4}, and claim that thereupon the conclusion of Theorem \ref{theo300} holds if we let
  \be{300.44}
    \Gamma_u(m,\eta):=\frac{m}{\pi} \cdot \frac{(b_0+\xi_0)^2}{b_0(b_0+\xi_0^2)}
  \ee
  and
  \be{300.45}
    \gamma(m,\eta):=\frac{m}{\pi} \cdot \frac{b_0}{b_0+\xi_0^2}
  \ee
  as well as
  \be{300.46}
    \Gamma_w(m,\eta):=2\Gamma_0.
  \ee
  To verify this, given $u_0$ and $w_0$ with the assumed properties we let $W_0$, $K_0$ and $U$ be defined by
  (\ref{W0}) and (\ref{31.1}), respectively, and fix $\UU$ as in (\ref{UU}), with $b(t):=b_0 e^{-\alpha t}$ and
  $a\in C^1([0,\infty))$ given by (\ref{a}).
  Then (\ref{W0}) and (\ref{300.04}) imply that
  \bas
    W_0(\xi)-\kz\xi
    &=& \int_0^{\sqrt{\xi}} \rho w_0(\rho)d\rho - \kz\xi \\
    &=& \bigg\{ \kz- \int_{\sqrt{\xi}}^1 \rho w_0(\rho)d\rho  \bigg\} - \kz\xi \\
    &=& (1-\xi)\kz - \frac{1-\xi}{2} \mint_{B_1\setminus B_{\sqrt{\xi}}} w_0 \\
    &\ge& (1-\xi)\kz - \frac{1-\xi}{2} \cdot \bigg\{ \mint_{B_1} w_0-\eta \bigg\} \\
    &=& (1-\xi) \cdot \bigg\{ \kz - \frac{1}{2} \mint_{B_1} w_0 + \frac{\eta}{2} \bigg\} \\
    &=& (1-\xi) \cdot \frac{\eta}{2}
    \qquad \mbox{for all } \xi \in (\xi_0,1),
  \eas
  because $\xi_0=R^2$ by (\ref{300.4}).
  Therefore,
  \be{300.66}
    \frac{W_0(\xi)-\kz\xi}{1-\xi}
    \ge \frac{\eta}{2}=\eta_0
    \qquad \mbox{for all } \xi \in (\xi_0,1),
  \ee
  so that Lemma \ref{lem2} applies to show that according to our choice of $\alpha$ and the fact that
  $b_0\in (0,\xi_0^2)$, taking $\parab$ as in (\ref{parab}) we have
  \be{300.7}
    \parab \UU(\xi,t) \le 0
    \qquad \mbox{for all $\xi\in (\xi_0,1)$ and } t>0.
  \ee
  We next combine (\ref{W0}) with (\ref{300.03}) and (\ref{300.46}) to see that
  \bea{300.77}
    \frac{W_0(\xi)}{\xi} - \kz
    &=& \frac{1}{\xi} \int_0^{\sqrt{\xi}} \rho w_0(\rho)d\rho
    -\int_0^1 \rho w_0(\rho)d\rho \nn\\
    &=& \frac{1}{2} \bigg\{ \mint_{B_{\sqrt{\xi}}} w_0 - \mint_{B_1} w_0 \bigg\} \nn\\
    &\ge& \frac{1}{2} \cdot \Gamma_w(m,\eta) \nn\\
    &=& \Gamma_0
    \qquad \mbox{for all } \xi\in (0,\xi_0),
  \eea
  again because $\xi_0=R^2$ by (\ref{300.4}).
  Consequently, Lemma \ref{lem4} asserts that
  \be{300.8}
    \parab \UU(\xi,t) \le 0
    \qquad \mbox{for all $\xi\in (0,\xi_0)$ and } t\in (0,t_0).
  \ee
  Moreover, since (\ref{300.66}) together with (\ref{300.77}) clearly implies that $W_0(\xi)-\kz\xi \ge 0$ for all
  $\xi\in (0,1)$, thanks to our choice of $b_0$ and $t_0$ and the fact that $\alpha<\alpha_\star$ we may employ
  Lemma \ref{lem1} to infer that
  \be{300.9}
    \parab \UU(\xi,t) \le 0
    \qquad \mbox{for all $\xi\in (0,\xi_0)$ and } t\ge t_0.
  \ee
  Now from the definition of $\UU$ and Lemma \ref{lem31} it is clear that
  \be{300.99}
    \UU(0,t)=U(0,t)=0
    \quad \mbox{and} \quad
    \UU(1,t)=U(1,t)=\frac{m}{2\pi}
    \qquad \mbox{for all } t>0.
  \ee
  In order to show that furthermore
  \be{300.10}
    \UU(\xi,0) \le U(\xi,0)
    \qquad \mbox{for all } \xi\in (0,1),
  \ee
  we first consider small values of $\xi$, for which from (\ref{300.01}) and (\ref{300.4}) we gain the inequality
  \bas
    U(\xi,0)
    = \int_0^{\sqrt{\xi}} \rho u_0(\rho)d\rho
    = \frac{\xi}{2} \mint_{B_{\sqrt{\xi}}} u_0
    \ge \frac{\xi}{2} \cdot \Gamma_u(m,\eta)
    \qquad \mbox{for all } \xi \in (0,\xi_0).
  \eas
  On the other hand, (\ref{UU}), (\ref{a}) and (\ref{300.44}) show that
  \bas
    \UU(\xi,0) = \frac{a(0)\xi}{b_0+\xi}
    \le \frac{a(0)\xi}{b_0}
    \le \frac{1}{2}\Gamma_u(m,\eta) \cdot \xi
    \qquad \mbox{for all } \xi \in (0,\xi_0),
  \eas
  and that hence (\ref{300.10}) is valid for any such $\xi$.
  For large $\xi$, by (\ref{a}) we have
  \bas
    \UU(\xi,0)
    &=& \frac{a(0) \cdot (b_0\xiü+\xi_0^2)}{(b_0+\xi_0^2)^2} \\
    &=& \frac{m}{2\pi} \cdot \frac{b_0\xi+\xi_0^2}{b_0+\xi_0^2} \\
    &=& \frac{m}{2\pi} \cdot \bigg\{ 1 - \frac{b_0}{b_0+\xi_0^2} \cdot (1-\xi)\bigg\}
    \qquad \mbox{for all } \xi \in (\xi_0,1),
  \eas
  whereas (\ref{300.02}), (\ref{300.45}) and again (\ref{300.4}) yield
  \bas
    U(\xi,0)
    &=& \int_0^1 \rho u_0(\rho)d\rho - \int_{\sqrt{\xi}}^1 \rho u_0(\rho)d\rho \\
    &=& \frac{m}{2\pi} - \frac{1-\xi}{2} \cdot \mint_{B_1\setminus B_{\sqrt{\xi}}} u_0 \\
    &\ge& \frac{m}{2\pi} - \frac{1-\xi}{2} \cdot \gamma(m,\eta) \\
    &=& \frac{m}{2\pi} - \frac{m}{2\pi} \cdot \frac{b_0}{b_0+\xi_0^2} \cdot (1-\xi)
    \qquad \mbox{for all } \xi \in (\xi_0,1).
  \eas
  We thereby conclude that the claimed ordering in (\ref{300.10}) indeed holds, so that on the basis of
  (\ref{300.7}), (\ref{300.8}), (\ref{300.9}) and (\ref{300.99}) we may invoke the comparison principle
  in Lemma \ref{lem12} to infer that $\UU(\xi,t)\le U(\xi,t)$ for all $\xi\in [0,1]$ and $t\ge 0$.
  In particular, this entails that for each fixed $t>0$ we must have
  \bas
    U_\xi(0,t) = \lim_{\xi\searrow 0} \frac{U(\xi,t)}{\xi}
    \ge \lim_{\xi\searrow 0} \frac{\UU(\xi,t)}{\xi}
    = \lim_{\xi\searrow 0} \frac{a(t)}{b(t)+\xi}
    =\frac{a(t)}{b(t)}
    = \frac{m}{2\pi} \cdot \frac{(b(t)+\xi_0)^2}{b(t)+\xi_0^2} \cdot \frac{1}{b(t)}.
  \eas
  Since $b(t)=b_0 e^{-\alpha t} \le b_0 < \xi_0^2$ for all $t>0$, this ensures that
  \bas
    u(0,t) = 2U_\xi(0,t)
    \ge \frac{m}{\pi} \cdot \frac{\xi_0^2}{2\xi_0^2} \cdot \frac{1}{b(t)}
    = \frac{m}{2\pi b_0} e^{\alpha t}
    \qquad \mbox{for all } t>0
  \eas
  and hence completes the proof.
\qed
\vspace*{10mm}
{\bf Acknowledgment.}
 The authors would like to thank the anonymous reviewers for their fruitful remarks.
 The first author is supported by the National Natural Science Foundation of China (No.
 11171061).
  This work was finished while the second author visited Dong Hua University in October 2013. He is grateful
  for the warm hospitality.
\end{document}